\documentclass[preprint,3p,12pt,numbers]{elsarticle}

\usepackage{amsmath}
\usepackage{amsfonts}
\usepackage{amscd}
\usepackage{amssymb}
\usepackage{mathrsfs}
\usepackage{amsthm}
\usepackage{ulem}
\usepackage{color}
\usepackage{soul}
\usepackage{needspace}
\usepackage{transparent}
\usepackage{enumitem}
\usepackage{natbib}
\usepackage[hidelinks]{hyperref}
\hypersetup{pdfauthor=author}

\newtheorem{Lemma}{Lemma}[section]
\newtheorem{Theorem}{Theorem}[section]
\newtheorem{Proposition}[Lemma]{Proposition}

\newtheorem{Definition}[Lemma]{Definition}
\newtheorem*{remark}{Remark}
\newtheorem{LemmaApdx}{Lemma}[]

\newcommand{\R}{\mathbb{R}}
\newcommand{\C}{\mathbb{C}}

\newcommand{\Z}{\mathbb{Z}}

\newcommand{\mbf}[1]{\mathbf{#1}}
\newcommand{\tl}[1]{\tilde{#1}}
\newcommand{\lp}{\left}
\newcommand{\rp}{\right}
\newcommand{\ri}{\mathrm{i}}
\renewcommand{\Re}{\mathrm{Re}}
\renewcommand{\Im}{\mathrm{Im}}
\newcommand{\p}{\partial}
\newcommand{\diag}{\,\,\mathrm{diag}}
\newcommand{\rs}{\mathrm{s}}
\newcommand{\ru}{\mathrm{u}}

\journal{J. Math. Anal. Appl.}

\begin{document}

\begin{frontmatter}
\title{Exponential Dichotomies for Elliptic Equations on Multidimensional Domains}  
\author[1]{Margaret Beck}
\author[1]{Ryan Goh}
\author[1]{Alanna Haslam-Hyde\corref{2}}
\cortext[2]{Corresponding author}
\ead{alanna.haslam-hyde@gordon.edu}
\affiliation[1]{organization={Boston University, Department of Mathematics and Statistics},
	addressline={665 Commonwealth Ave}, 
	city={Boston},
	state={MA},
	postcode={02215}, 
	country={USA}}
\begin{abstract}
	The existence of exponential dichotomies has been well-established as a powerful tool to study existence, stability, and bifurcations of coherent structures. Currently, the application of exponential dichotomies to elliptic problems posed on multi-dimensional domains is predominately limited to the context of cylindrical spatial domains. Recent work in \cite{Beck21-SES} has shown how to extend the method of spatial dynamics, in which one views a spatial variable as a time-like evolutionary variable, to general multi-dimensional spatial domains. In this paper, we show that exponential dichotomies exist for a class of spatial dynamical systems arising in this more general setting, thus allowing for their use in future analyses of coherent structures.
\end{abstract}
\end{frontmatter}

\section{Introduction}
\setcounter{section}{1}
A current challenge in the study of pattern formation is the intractability of analytic dynamical systems approaches for pattern-forming partial differential equations in higher spatial dimensions. Dynamical systems techniques, for example center manifold reductions, geometric singular perturbation theory, and exponential dichotomies, have been applied to study the existence, stability, and bifurcations of coherent structures in a variety of settings posed on one-dimensional spatial domains, but their use in higher dimensions still remains limited to certain contexts. It is the aim of this work to broaden the class of higher-dimensional settings where such analyses can take place by furthering the theory of exponential dichotomies.

\vspace{6pt}\noindent
\textbf{Spatial dynamics for one-dimensional spatial domains.}
To study elliptic equations from a dynamical systems perspective, we employ the technique of spatial dynamics. In the simplest context, spatial dynamics is the trick of rewriting a second order differential equation into a system of first order equations. For example, the equation $u_{xx} + F(u) = 0 $, where $u(x)\in \C^d$ for each fixed $x\in\R$, becomes the system
\begin{align*}
	u_x & = v \\
	v_x & = - F(u).
\end{align*}
where $(u(x),v(x))\in\C^{d}\times \C^d$. 
This same procedure provides a method for studying the stationary and traveling-wave problems associated with PDEs posed on $\R$. For example, a traveling wave ansatz $u(\xi)=u(x+ct)$ applied to the reaction-diffusion equation $u_t=u_{xx}+F(u)$ yields $c u_{\xi} = u_{\xi\xi}+F(u)$ which can be written as a first order system of ODEs as
\begin{equation}\label{eq.traveling_wave}
	\begin{split}
		u_\xi &= v\\
		v_\xi &= cv - F(u),
	\end{split}
\end{equation}
where $(u(\xi),v(\xi))\in\C^{d}\times \C^d$.

Equation \eqref{eq.traveling_wave} is the classical setting in which ODE tools can be explicitly applied to study pattern formation. For example, bifurcation theory, center manifold reduction, and geometric singular perturbation theory are traditional ODE tools which have been extensively applied in this setting to understand the emergence of non-uniform stationary spatial patterns on $\R$; see \cite{Doelman19} for a helpful introduction. 
Of particular interest in this paper is the tool of exponential dichotomies, which were originally developed to extend the notion of invariant stable and unstable linear subspaces to the non-autonomous setting \cite{Coppel78}. Because they characterize forward and backward exponential decay of linear (non-autonomous) systems, exponential dichotomies can play a significant role in constructing and analyzing bounded solutions that decay to known backgrounds states in the tails. Notably, exponential dichotomies can be used to construct bounded perturbations of pulses and fronts, which appear as homoclinic and heteroclinic solutions of \eqref{eq.traveling_wave}, as in \cite{Palmer84}. Other work on exponential dichotomies has shown their connection with the Fredholm properties of related differential operators; they therefore also play a significant role in understanding the temporal stability properties of stationary and traveling wave problems (e.g. \cite{Sandstede02}).

\vspace{6pt}\noindent
\textbf{Spatial dynamics for cylindrical domains.}
Motivated by the success of dynamical systems tools for studying PDEs with one spatial dimension,
\cite{Kirchgassner82} extended the notion of spatial dynamics to PDEs posed on cylindrical spatial domains of the form $\R\times\Omega\subset\R^n$ where $\Omega$ is a compact subset of $\R^{n-1}$. In this channel setting, stationary or traveling-wave solutions can be studied as solutions to a spatial evolutionary equation which describes how functions on the compact cross-section $\Omega$ evolve in the unbounded variable $x\in\R$. Consequently, spatial dynamics on the channel becomes an infinite-dimensional evolutionary equation rather than a finite-dimensional ODE in the case of one-dimensional spatial domains. For example, stationary solutions of $u_t = \Delta u + F(u)$ posed on a general $n$-dimensional channel $\R\times \Omega \subset \R^n$, will satisfy
\begin{equation}\label{eq.channel_pde}
	  u_{xx} + \Delta_{\Omega}u + F(u)=0.
\end{equation}
By writing this equation as a first order system in the unbounded variable $x$, we can study the stationary problem as an evolutionary equation
\begin{equation} \label{eq.channel}
	\begin{split}
		u_x &= v\\
		v_x &= -\Delta_{\Omega}u - F(u)
	\end{split}
\end{equation}
posed on a Banach space, so that for each fixed $x\in \R$, $(u(x),v(x)) \in H^1(\Omega,\C^d)\times L^2(\Omega,\C^d)$.

This approach yields linearized systems with compact resolvent, allowing for the application of center manifold techniques to reduce the infinite-dimensional dynamics of \eqref{eq.channel} to a finite-dimensional manifold governing small bounded solutions in both autonomous \cite{Kirchgassner82} and non-autonomous systems \cite{Mielke86}.
A similar finite-dimensional reduction by \cite{Mielke94} extended these techniques to encompass not just small amplitude solutions but also larger bounded solutions. Front interfaces between these larger bounded solutions were also studied by \cite{Mielke94} via an extension of the heteroclinic front construction techniques used in the analogous finite-dimensional setting.

Further extension of classic dynamical systems techniques to spatial dynamics on the channel proved challenging due to the ill-posedness of \eqref{eq.channel}. The existence of well-behaved solutions in either forward or backward evolution is not guaranteed for these spatial evolutionary equations since the linear operator $\left(\begin{smallmatrix}
	0 & 1\\ -\Delta_\Omega &0 
\end{smallmatrix}\right)$ has an unbounded point spectrum with real part going to both positive and negative infinity.
The theory of exponential dichotomies was extended to this setting by \cite{PeterhofSandstedeScheel97}, and is a primary tool to overcome the challenge of ill-posedness. 
Just as in the finite dimensional setting, exponential dichotomies describe how to split the underlying space into invariant families of stable and unstable subspaces. Consequently, even for ill-posed evolutionary equations, exponential dichotomies describe where well-behaved solutions exist and decay exponentially in the forward direction (the stable family of subspaces) and where solutions exist and decay exponentially in the backward direction (the unstable family of subspaces). 
Furthermore, \cite{PeterhofSandstedeScheel97} showed how exponential dichotomies for infinite-dimensional problems can be combined with Lyapunov-Schmidt techniques in the style of \cite{Palmer84} to construct bounded solutions of nonlinear equations that remain close to a pulse solution. 

The primary reason why spatial dynamics works for studying elliptic PDEs on cylindrical domains is because the resulting linearized operator has compact resolvent. In the setting of \eqref{eq.channel} this is due to the fact that the cross-sectional spatial domain $\Omega$ is a compact subset of $\R^{n-1}$, but it can also be achieved in other ways. For example, a linearized operator with compact resolvent is obtained when studying time-periodic solutions of a reaction-diffusion equation posed on $\R$ by considering how functions in $L^2_{\mathrm{per}}([0,T],\C^d)$ evolve in the unbounded spatial variable $x\in \R$ (e.g.  \cite{SandstedeScheel01,SandstedeScheel01_structurespectra,SandstedeScheel04_classify}).
A similar effect is also achieved when studying traveling fronts between spatially periodic background states. In this case, by treating the direction of propagation as the evolutionary variable, it is the imposition of periodicity in the cross-sectional space that ensures the linearized operator has compact resolvent. For example, \cite{EckmannWayne91} studied one-dimensional traveling fronts in which periodic patterns invade the trivial state, and \cite{Doelman03} studied planar traveling fonts in which stable hexagonal patterns invade various unstable patterns.
	
\noindent
\textbf{Spatial dynamics on more general multidimensional domains.}
Spatial dynamics was extended to more general two-dimensional patterns, including target \cite{Scheel03} and spiral waves \cite{Scheel98}, by decomposing $\R^2$ into spherical coordinates $(r,\varphi)\in\R^+\times S^{1}$ with the radial variable $r$ playing the role of the unbounded evolutionary variable.
For example, radially symmetric stationary solutions of $u_t=\Delta u + F(u)$ satisfy $\Delta u + F(u) = 0$, which when written in spherical coordinates (and imposing $u(r,\varphi)\equiv u(r)$) becomes
\begin{equation}\label{eq.radial_pde}
\quad u_{rr} + \left(\frac{n-1}{r}\right)u_r + F(u)=0.
\end{equation}
This can then be written as a first order system in the variable $r$ as
\begin{align*}
	u_r &= v\\
	v_r &= -\left(\frac{n-1}{r}\right)v - F(u),
\end{align*}
where for each fixed $r\in\R^+$, $(u(r),v(r))\in\C^d\times \C^d$. The results of \cite{Scheel03} overcome both the non-autonomous dependence on $r$ and the singularity at $r=0$, developing a systematic bifurcation theory including extensions of center manifold reduction and normal form theory.

A subsequent study of spiral waves from a spatial dynamics perspective was given by \cite{SandstedeScheel21} in which the existence of exponential dichotomies were proven and applied to rigorously understand the spectral properties and robustness of spiral waves on both bounded and unbounded domains. In addition, \cite{SandstedeScheel21} used this perspective to characterize the behavior of eigenfunctions to the linearization which indicate the spatial structure of patterns bifurcating from radially symmetric spiral waves.

Recent work, \cite{Beck21-SES}, has extended the use of spatial dynamics further to arbitrary Euclidean spatial domains including, of course, $\R^n$. The application of spatial dynamics in $\R^n$ looks much the same as in the setting of $\R^2$, relying again on spherical coordinates $(r,\varphi)\in \R^+\times S^{n-1}$. The elliptic equation $0 = \Delta u - V(x)u$ with $V:\R^n\to \C^{d\times d}$ and $u(x)\in\C^d$ for $x\in\R^n$ is rewritten in spherical coordinates as 
\begin{equation}\label{eq.spherical_pde}
u_{rr}+\left(\frac{n-1}{r}\right)u_r + \frac{1}{r^2}\Delta_{S^{n-1}} u  - V(r,\varphi)u=0, \qquad u(r,\cdot)\in L^2(S^{n-1},\C^d)
\end{equation}
where $\Delta_{S^{n-1}}$ is the spherical Laplacian (well understood via spherical harmonics). The above can be expressed as a first order system in $r$ as 
\begin{equation}
	\mathbf{u}_r = \mathcal{A}(r)\mathbf{u} 
	= \begin{pmatrix} 0 & 1 \\ V(r,\varphi) - \frac{1}{r^2} \Delta_{S^{n-1}} & - \frac{(n-1)}{r} \end{pmatrix} \mathbf{u}
	\label{eq.spatial_eq}
\end{equation}
posed on the Banach space $X$ defined as
\begin{equation}\label{eq.X_norm}
	X = H^1(S^{n-1},\C^d)\times L^2(S^{n-1},\C^d), \qquad 	\big\|\mathbf{u}(r,\cdot)\big\|_{X}^2 = \big\|u(r,\cdot)\big\|_{H^1}^2 + \big\|v(r,\cdot)\big\|_{L^2}^2
\end{equation}
so that, for each $r\in\R^+$, $\mbf{u}(r) = (u(r),v(r))\in X$. 

The results of \cite{Beck21-SES} establish not only the method of spatial dynamics, but also the equivalence of weak PDE solutions and solutions to the corresponding spatial evolutionary equation in both $\R^n$ and general bounded Euclidean domains.
These results lay the groundwork for adapting classical dynamical systems techniques to this newer spatial dynamics setting. In fact, \cite{Beck21-ExpDich} showed the existence of an exponential dichotomy for $r\in (0,R]$ (the core field), overcoming the singularity at $r=0$ in a relatively straightforward way using the reparameterization $\tau=\ln r$ and an exponential rescaling of the resulting system posed on $\tau\in \R^-$ (see also \cite{Scheel98}).

\begin{remark}
	Though this discussion has predominately been framed in terms of the existence of solutions, questions of stability are also often addressed using the techniques described above. For example, one could study an eigenvalue problem for the linearization about a coherent structure in any of these settings by adding the linear term $-\lambda u$ to Equations \eqref{eq.channel_pde}, \eqref{eq.radial_pde}, or \eqref{eq.spherical_pde}.
\end{remark}

\noindent
\textbf{Our contribution: Far-field exponential dichotomies for elliptic problems posed on $\R^n$.}
To extend the results of \cite{Beck21-ExpDich} to problems posed on the unbounded domain $\R^n$, we prove the existence of a far-field dichotomy---an exponential dichotomy on the interval $r\in[1,\infty)$---for the same spatial evolutionary equation \eqref{eq.spatial_eq}. The basis of our proof relies on the robustness of exponential dichotomies, the property that exponential dichotomies persist through sufficiently small perturbations. 
While the standard robustness result in the cylindrical setting given by \cite{PeterhofSandstedeScheel97} can be applied in a relatively direct way for the core-field dichotomy proof in \cite{Beck21-ExpDich}, 
our proof of the far-field dichotomy is complicated by the non-uniform unboundedness of the $r^{-2}\Delta_{S^{n-1}}$ term in \eqref{eq.spatial_eq}. 

Though the $r^{-2}\Delta_{S^{n-1}}$ term vanishes as $r$ is sent to infinity and does not primarily control the hyperbolicity of the system, its unboundedness with respect to the $L^2$ norm means that it cannot simply be considered a perturbation. Instead, we must show that the principal component of \eqref{eq.spatial_eq}, including this spherical Laplacian term, has an exponential dichotomy in an appropriately chosen function space.
Taking inspiration from the robustness argument in \cite{Scheel98} and \cite{SandstedeScheel21}, we choose to work in a function space with an explicitly $r$-dependent norm which complements the non-uniformity of the $r^{-2}\Delta_{S^{n-1}}$ term by ensuring that it remains bounded while not dominating the asymptotic hyperbolicity of the system. 

In particular, observing that $\mathcal{A}(r)$ is both closed and densely defined in the space $X$ with domain $X^1 = H^2(S^{n-1})\times H^1(S^{n-1})$, we adapt the natural $X$ norm, as defined in \eqref{eq.X_norm}, to an $r$-dependent norm by weighting the $H^1$ component:
\begin{align}
	\begin{array}{c}
		X_r = H^1(S^{n-1},\C^d)\times L^2(S^{n-1},\C^d), \\[1em]
		\big\|\mathbf{u}(r,\cdot)\big\|_{X_r}^2 := \frac{1}{r^2}\big\|u(r,\cdot)\big\|_{H^1}^2 + \big\|u(r,\cdot)\big\|_{L^2}^2 + \big\|v(r,\cdot)\big\|_{L^2}^2
	\end{array}\label{eq.Xr_norm}
\end{align}
In this $r$-dependent function space $X_r$, the dynamics of the asymptotic hyperbolic operator dominate norm estimates as $r$ goes to infinity while the influence of the spherical Laplacian remains bounded under the weighted $H^1$ norm. This choice of norm turns out to be the right function space in which to show our system has an exponential dichotomy. A carefully constructed, norm-preserving transformation between $X_r$ and $L^2(S^{n-1})\times L^2(S^{n-1})$---presented in Lemma \ref{l.isomorphism}---is the key trick to applying robustness in this setting and allows us to prove our result.

It is worth noting that both \cite{Scheel98} and \cite{SandstedeScheel21} design similar $r$-weighted norms 
to establish the existence of exponential dichotomies in the spiral wave settings.
 Our choice of $r$-dependent norm differs from theirs in the middle unweighted $L^2$ norm of $u$ where they instead have an unweighted $H^{1/2}$ norm. This is due to an angular advection term that appears in the spiral waves context.  
Beyond the absence of an advection term in our setting, the main distinction between the proof of exponential dichotomies in the spiral waves setting and the proof presented here is the way that we decompose the system into a principal part and a perturbation. 
In \cite{SandstedeScheel21}, the principal component includes only the differential operators, treating the additional terms which are analogous to our asymptotic potential, $V_\infty = \lim_{r\to\infty}V(r,\varphi)$, as perturbations. 
The principal component in our analysis includes both the spherical Laplacian term and the asymptotic potential.
As a result, our proof involves slightly more technical bounds to prove Lemma \ref{l.isomorphism} (comparable to Lemma 5.6 in \cite{SandstedeScheel21}), while the proof in \cite{SandstedeScheel21} requires a more elaborate argument to bound the perturbation terms.

\vspace{6pt}\noindent
\textbf{Organizational overview of paper.}
In section \ref{sec:2}, we present this framework for the existence of a far-field dichotomy for equation \eqref{eq.spatial_eq} in greater detail. Later, in Section \ref{sec:3}, we prove the Lemmas and existence results stated in Section \ref{sec:2}, with many of the technical details contained in the Appendix. Finally, in Section \ref{sec:4}, we discuss possible applications of our results, particularly with respect to the existence and construction of localized spatial patterns in $\R^n$.

\section{Definitions and Results.}\label{sec:2}
\setcounter{section}{2}

Our approach to proving the existence of exponential dichotomies relies heavily on the property of robustness, the persistence of exponential dichotomies through sufficiently small perturbations. Therefore, to prove that \eqref{eq.spatial_eq} has an exponential dichotomy on an appropriately chosen function space, we decompose it into a principal part, $\tilde{\mathcal{A}}(r)$, that we can directly prove has an exponential dichotomy, and a perturbation, $\tilde{B}(r)$, that we can bound appropriately. This decomposition can be characterized as
\begin{equation*}
	\mathcal{A}(r) = 
	\underbrace{
		\underbrace{\begin{pmatrix}0 & 1\\ V_\infty & 0 \end{pmatrix}}_{\mathcal{A}_\infty}
		+ \begin{pmatrix}	0 & 0 \\ -\frac{1}{r^2}\Delta_{S^{n-1}} & 0	\end{pmatrix}
	}_{\tilde{\mathcal{A}}(r)} 
	+ \underbrace{
		\vphantom{
			\underbrace{\begin{pmatrix} 0 & 0 \\ 0 & 0 \end{pmatrix}}_{\mathcal{A}_\infty}}
		\begin{pmatrix}	0 & 0 \\ V(r,\cdot)-V_\infty &-\frac{(n-1)}{r}	\end{pmatrix}
	}_{\tilde{B}(r)}
\end{equation*}
where $\mathcal{A}_\infty = \lim_{r\to\infty}\mathcal{A}(r)$ is the asymptotic operator of $\mathcal{A}(r)$ under the assumption that $V(r,\varphi)$ decays to a constant potential $V_\infty \in \C^{d\times d}$ as $r$ goes to infinity.

Often the asymptotic operator is a suitable choice for the principal part of a robustness argument. Since spectral hyperbolicity is equivalent to the existence of exponential dichotomies in a constant coefficient system, it is straightforward to look for an exponential dichotomy of the asymptotic system:
\begin{equation}
	\mbf{u}_r = \begin{pmatrix}0 & 1\\ V_\infty & 0 \end{pmatrix}\mbf{u} = \mathcal{A}_\infty \mbf{u}. \label{eq.asymp_eq}
\end{equation}	
In this context, however, it is not sufficient to consider $\mathcal{A}_\infty$ as the principal part, as we cannot directly apply a robustness argument to show that \eqref{eq.spatial_eq} has an exponential dichotomy when \eqref{eq.asymp_eq} does. In particular,
if we were to take $\mathcal{A}_\infty$ as the principal operator, we would be left with a remainder that is unbounded relative to $\mathcal{A}_\infty$ under the natural product norm on $X$ due to the spherical Laplacian term. Since the spherical Laplacian term can therefore not be left as a perturbation, we must include it in our principal operator and instead choose a norm that both bounds the Laplacian when $r$ is small and highlights the effects of $V_\infty$ when $r$ is large. This is accomplished by the norms defined in \eqref{eq.Xr_norm}.

With this new weighted family of function spaces, $X_r$ in \eqref{eq.Xr_norm}, we can prove the existence of an exponential dichotomy for a principal part containing 
the leading order derivatives. This principal part is given by
\begin{equation}\label{eq.principal_op}
	\tilde{\mathcal{A}}(r) = \begin{pmatrix}
		0 & 1\\ V_\infty - \frac{1}{r^2}\Delta_{S^{n-1}} & 0
	\end{pmatrix}
\end{equation}
and leaves a remaining perturbation given by
\begin{equation}\label{eq.perturb_op}
	\tilde{B}(r) = \begin{pmatrix}
		0 & 0\\ V(r,\varphi)-V_\infty & -\frac{(n-1)}{r}
	\end{pmatrix}.
\end{equation}
We will later refer to the principal equation:
\begin{equation}
	\mbf{u}_r = \tilde{\mathcal{A}}(r) \mbf{u}. \label{eq.principal_eq}
\end{equation}

Before further discussing the existence proof, we must first define what is meant by an exponential dichotomy on the $r$-dependent space $X_r$. 

	\begin{remark}
		The choice of an $r$-dependent norm effectively creates a family of function spaces $X_r$ where the subscript denotes which value of $r$ is to be used in the norm. Consequently, a function $u : S^{n-1}\to\C^d$ may live in any $X_r$, however we will typically consider $u : [0,\infty)\times S^{n-1}\to \C^d$ evaluated at a fixed value of $r$ as living in the corresponding $X_r$ space; for example, $u(r,\cdot)\in X_r$ and $u(\rho,\cdot)\in X_\rho$. For this reason, we refer to $X_r$ both as a family of function spaces and as a single function space.
	\end{remark}

 \begin{Definition}\label{d.ed-Xr} \textup{(Definition 5.2 in \cite{SandstedeScheel21})}
 	We say that a radial spatial evolutionary equation of the form $\mbf{u}_r = \mathcal{L}(r)\mbf{u}$ has an exponential dichotomy in the $r$-dependent family of Banach spaces $X_r$ on a sub-interval $J\subseteq \R^+$ if there exist two strongly continuous families of projections $P^{\mathrm{s},\mathrm{u}}(r)\in L(X_{r})$ with $P^\mathrm{s}(r)+P^\mathrm{u}(r) = I$, two families of linear operators $\Phi^{\mathrm{s},\mathrm{u}}(r,\rho)\in L(X_\rho, X_r)$, and positive constants $C$ and $\eta$ such that the following statements are true: 
 \begin{itemize}
 	\item (Stability)
 	For all $\rho\in J$ and $\mathbf{u}_0\in X_\rho$, there exists a solution $\Phi^{\rs}(r,\rho) \mathbf{u}_0$, defined and continuous for both $r$ and $\rho$ on $r\ge\rho$, and differentiable in both $r$ and $\rho$ for $r>\rho$. The operator $P^{\rs}(\rho)\mathbf{u}_0 := \Phi^{\rs}(\rho,\rho) \mathbf{u}_0$ gives the stable projection, and
 	\begin{equation*}\label{eq.stablity_Xr}
 		\|\Phi^\rs(r,\rho) \mathbf{u}_0\|_{X_r}\leq C e^{-\eta(r - \rho)} \|\mathbf{u}_0\|_{X_\rho}, \qquad r\geq \rho, \,\, r,\rho\in J,
 	\end{equation*}
 	where $C$ and $\eta$ are independent of $r$, $\rho$, and $\mathbf{u}_0$.
 	\item (Instability) 
 	For all $\rho\in J$ and $\mathbf{u}_0\in X_\rho$, there exists a solution $\Phi^{\ru}(r,\rho) \mathbf{u}_0$, defined and continuous for both $r$ and $\rho$ on $r\le\rho$, and differentiable in both $r$ and $\rho$ for $r<\rho$. The operator $P^{\ru}(\rho)\mathbf{u}_0 := \Phi^{\ru}(\rho,\rho) \mathbf{u}_0$ gives the unstable projection, and 
 	\begin{equation*}\label{eq.instability_Xr}
 		\|\Phi^\ru(r,\rho) \mathbf{u}_0\|_{X_r}\leq C e^{\eta(r - \rho)} \|\mathbf{u}_0\|_{X_\rho}, \qquad r\leq \rho, \,\, r,\rho\in J,
 	\end{equation*}
 	where $C$ and $\eta$ are independent of $r$, $\rho$, and $\mathbf{u}_0$. 
 	\item (Invariance) The solutions $\Phi^{\rs}(r,\rho)\mathbf{u}_0$ and $\Phi^\ru(r,\rho)\mbf{u}_0$ are invariant under the respective projections. That is, $\Phi^{\rs}(r,\rho)\mathbf{u}_0\in \mathrm{Rg}\, P^{\rs}(r)$ for all $r\geq \rho\in J$ and $\Phi^{\ru}(r,\rho)\mathbf{u}_0\in \mathrm{Rg}\, P^{\ru}(r)$ for all $r\leq \rho\in J$. Equivalently,
 	\begin{align*}
 		&\,\, P^\rs(r)\Phi^\rs(r,\rho) = \Phi^\rs(r,\rho) = \Phi^\rs(r,\rho)P^\rs(\rho)\\
 		\quad &P^\ru(r)\Phi^\ru(r,\rho) = \Phi^\ru(r,\rho) = \Phi^\ru(r,\rho)P^\ru(\rho).
 	\end{align*}
 \end{itemize}
 \end{Definition}

Throughout this paper, we also consider exponential dichotomies posed on the fixed space $L^2(S^{n-1},\C^d)\times L^2(S^{n-1},\C^d)$; we omit the standard definition in the context of a fixed Banach space, but the interested reader can find one in \cite{PeterhofSandstedeScheel97}, noting that it is analogous to the definition above.

For completeness, we now also state the robustness of exponential dichotomies in the context of the $r$-dependent family of function spaces $X_r$.

\begin{Proposition}\textup{(Proposition 5.9 in \cite{SandstedeScheel21})}\label{p.roughness}
	Exponential dichotomies are robust under sufficiently small perturbations. More precisely, if a system of the form $\mathbf{u}_r = \mathcal{L}(r) \mathbf{u}$ has an exponential dichotomy on some interval $J\subseteq \R^+$ in the $r$-dependent family of $X_r$ spaces with constants $C,\eta >0$, then for all $\varepsilon>0$ and $\tl\eta\in(0,\eta)$ there are constants $\tl C, \delta>0$ such that the perturbed system 
	$$\mathbf{u}_r = (\mathcal{L}(r)+ \mathcal{B}(r)) \mathbf{u}$$
	with
	$$\| \mathcal{B}(r)\|_{L(X_r)} \leq \delta, \qquad r\in J$$
	has an exponential dichotomy on $J$ in the same family of spaces, $X_r$, with constants $\tl C$ and $\tl \eta$ and with projections $\varepsilon$-close to the projections of $\mathbf{u}_r = \mathcal{L}(r) \mathbf{u}$. Furthermore, if $\|\mathcal{B}(r)\|_{L(X_r)}  = O(1/r)$, then one can take $\tl\eta = \eta$.
\end{Proposition}

The comparable standard robustness result in the setting of a fixed Banach space is given as Theorem 1 in \cite{PeterhofSandstedeScheel97}. Throughout the paper, we will invoke robustness in both settings.

To apply Proposition \ref{p.roughness}, we must first establish the existence of an exponential dichotomy for the unperturbed or principal system, in this case \eqref{eq.principal_eq}. 
Making reasonable assumptions on $V_\infty$, the asymptotic system \eqref{eq.asymp_eq} will have an exponential dichotomy on $L^2(S^{n-1},\C^d)\times L^2(S^{n-1},\C^d)$. The main difficulty lies in applying robustness to extend the exponential dichotomy for \eqref{eq.asymp_eq} to one for \eqref{eq.principal_eq} on the family of $X_r$ spaces.
The majority of the proof is therefore dedicated to this first application of robustness. Once we have an exponential dichotomy for the principal system \eqref{eq.principal_eq}, we then apply robustness a second time to obtain an exponential dichotomy on the full system \eqref{eq.spatial_eq}.

The second application of robustness---to prove an exponential dichotomy for \eqref{eq.spatial_eq} given an exponential dichotomy for \eqref{eq.principal_eq}---is direct, assuming that $\tilde{B}(r)$ is sufficiently small, and in fact, decaying like $\mathcal{O}(1/r)$. This will hold if the potential function $V(r,\varphi)$ decays like $\mathcal{O}(1/r)$ to the asymptotic potential $V_\infty$. We therefore make the following hypothesis:
\begin{enumerate}[leftmargin=0.7in,rightmargin=0.2in]
	\item[\textbf{(H1.)}] 
	There exists a constant $C>0$ such that for $r\ge 1$, $\|V(r,\varphi)-V_\infty\|\le C/r$.
\end{enumerate}
Turning to the existence of an asymptotic exponential dichotomy, we note that since $\mathcal{A}_\infty$ is autonomous and homogeneous across $S^{n-1}$, 
we can explicitly construct an exponential dichotomy for \eqref{eq.asymp_eq} via a spectral decomposition whenever the point spectrum of $\mathcal{A}_\infty$ does not intersect the imaginary axis $\ri \R\subset\C$. This will hold whenever the eigenvalues of $V_\infty$ do not lie on the negative real line, $(-\infty,0]$, so we therefore make the following hypothesis:
\begin{enumerate}[leftmargin=0.7in,rightmargin=0.2in]
	\item[\textbf{(H2.)}] 
	For each eigenvalue $\lambda_\ell$ of $V_\infty$, $\lambda_\ell \in \C\backslash(-\infty,0]$.
\end{enumerate}
Finally, we also assume that $V_\infty$ has only semi-simple eigenvalues, and is therefore diagonalizable. This last hypothesis greatly simplifies the proofs to follow, allowing for an explicit spectral decomposition of $\tilde{\mathcal{A}}(r)$ into scalar pieces. 
In \S \ref{sec:spherical_harmonics} we will carry out this decomposition explicitly and comment on the extendability of related arguments to the case of nontrivial Jordan blocks.
\begin{enumerate}[leftmargin=0.7in,rightmargin=0.2in]
	\item[\textbf{(H3.)}] $V_\infty$ is diagonalizable with eigenvalues $\{\lambda_\ell\}_{\ell=1}^d$ and corresponding eigenvectors $\{\rho_\ell\}_{\ell=1}^d \subset \C^d$.
\end{enumerate}

Hypotheses (H2) and (H3) allow us to fully characterize an exponential dichotomy for the asymptotic system \eqref{eq.asymp_eq}. 
We must then understand the relationship between the principal operator $\tilde{\mathcal{A}}(r)$, in \eqref{eq.principal_op}, and the asymptotic operator $\mathcal{A}_\infty$. To do this, we establish an isomorphism, or rather a family of isomorphisms, that can translate between the $X_r$ setting and the $L^2(S^{n-1})\times  L^2(S^{n-1})$ setting. The basis of this family of isomorphisms, established in Lemma \ref{l.isomorphism}, is the differential operator defined by
\begin{equation}
	A(r):= V_\infty - \frac{1}{r^2}\Delta_{S^{n-1}}. \label{eq.Ar_iso}
\end{equation}
Note that with this new notation, the principal operator \eqref{eq.principal_op} can be written as $\tilde{\mathcal{A}}(r) = \left(\begin{smallmatrix}
	0 & 1\\ A(r) & 0
\end{smallmatrix}\right)$.

\begin{Lemma}\label{l.isomorphism}
	If \textup{(H2)} and \textup{(H3)} hold then the following statements are true:
	\begin{enumerate}[label=(\roman*)]
		\item The family of operators $A(r)^{1/2}:H^1(S^{n-1})\rightarrow L^2(S^{n-1})$, where $A(r)$ is given in \eqref{eq.Ar_iso}, is well-defined over all values of $r\ge1$. In addition, the family is uniformly bounded, strongly continuous, and continuously differentiable in $r$ for $r>1$, with
		$$
		\left\| \p_r\!\lp[A(r)^{1/2}\rp]\!A(r)^{-1/2}\right\|_{L(L^2)} \leq C r^{-1}
		$$
		for some constant $C>0$ independent of $r$.
		\item For each $r\ge 1$, the operator $(u,v)\mapsto (A(r)^{1/2}u,v)$ defined on $X_r\rightarrow L^2(S^{n-1})\times L^2(S^{n-1})$ is an isomorphism and satisfies
		\begin{equation*}\label{e:Xr_norm_equiv}
			C_1 \lp( \|A(r)^{1/2} u\|_{L^2}^2 + \|v\|_{L^2}^2\rp)^{1/2}\leq \|\mathbf{u}\|_{X_r}\leq C_2\lp( \|A(r)^{1/2} u\|^2_{L^2} + \|v\|_{L^2}^2\rp)^{1/2}
		\end{equation*}
		for all $\mathbf{u} = (u,v) \in X_r$ with constants $C_1,C_2>0$ independent of both $r$ and $\mbf{u}$.
	\end{enumerate}
\end{Lemma}

Lemma \ref{l.isomorphism} establishes $(u,v)\mapsto(A(r)^{1/2}u,v)$ as an isomorphism from $X_r$ to $L^2 \times L^2$.
Applying the corresponding change of variables $\tilde{u} = A(r)^{1/2}u$ to the principal equation \eqref{eq.principal_eq} we get the following equation posed on $L^2(S^{n-1},\C^d)\times L^2 (S^{n-1},\C^d)$: \\
\begin{equation}\label{eq.COV_eq}
	\mbf{\tilde{u}}_r = \begin{pmatrix}
		\tilde{u}_r\\v_r
	\end{pmatrix} = \begin{pmatrix}
		A(r)^{1/2}v + \partial_r\!\lp[A(r)^{1/2}\rp]A(r)^{-1/2}\tilde{u}\\ A(r)^{1/2}\tilde{u}
	\end{pmatrix} = \begin{pmatrix}
		0& A(r)^{1/2}\\ A(r)^{1/2} & 0
	\end{pmatrix}\mbf{\tilde{u}} + \mathcal{O}\lp(1/r\rp)\mbf{\tilde{u}}
\end{equation}

By Proposition \ref{p.roughness} and Lemma \ref{l.isomorphism}, an exponential dichotomy for \eqref{eq.spatial_eq} in $X_r$ can be established by proving the existence of an exponential dichotomy for the principal part of \eqref{eq.COV_eq} in $L^2\times L^2$.
	
We proceed by performing another change of variables, $w^\pm = \tilde{u}\pm v$, to diagonalize the leading order part of \eqref{eq.COV_eq}:
	\begin{equation}\label{eq.diagonalized}
		\mbf{w}_r = \begin{pmatrix}
			w^+_r\\w^-_r
		\end{pmatrix} = \begin{pmatrix}
			A(r)^{1/2}w^+\\-A(r)^{1/2}w^-
		\end{pmatrix} = \begin{pmatrix}
			A(r)^{1/2} & 0 \\ 0 & -A(r)^{1/2}
		\end{pmatrix} \mbf{w}.
	\end{equation} 
	
	Once we arrive at \eqref{eq.diagonalized}, it is relatively straightforward to explicitly construct an exponential dichotomy using the spectral decomposition of the spherical Laplacian and the diagonalization of $V_\infty$.
	
	\begin{Lemma}\label{l.ED_diagonalized}
		Under hypotheses \textup{(H1)}, \textup{(H2)}, and \textup{(H3)}, \eqref{eq.diagonalized} has an exponential dichotomy,\, $\tilde{\Phi}^{\rs/\ru}(r,t)$, in $L^2(S^{n-1})\times L^2(S^{n-1})$ on the interval $r\ge1$ with the same decay rate as the asymptotic system, \eqref{eq.asymp_eq}: \,$\eta = \displaystyle{\tilde{\eta} = \min_{1\le\ell\le d}\mathrm{Re}[\lambda_\ell^{1/2}]}$.
	\end{Lemma}
	
	This provides the last stepping stone before our primary result:
	
	\begin{Theorem}\label{thrm.ED_existence}
		Under hypotheses \textup{(H1)}, \textup{(H2)}, and \textup{(H3)}, \eqref{eq.spatial_eq} has an exponential dichotomy in $X_r$ for $r\ge1$ with the same decay rate as the asymptotic equation \eqref{eq.asymp_eq}: $\displaystyle{\eta = \tilde{\eta} = \min_{1\le \ell\le d}\mathrm{Re}[\lambda_\ell^{1/2}]}$.
	\end{Theorem}
	
	Lemma \ref{l.ED_diagonalized} will be proven in \S \ref{sec:lem_proofs} by applying robustness of exponential dichotomies in $L^2\times L^2$ (Theorem 1 in \cite{PeterhofSandstedeScheel97}). Then, Theorem \ref{thrm.ED_existence} will be proven in \S \ref{sec:thrm_proof} by using Lemmas \ref{l.isomorphism} and \ref{l.ED_diagonalized} to apply robustness of exponential dichotomies in $X_r$ (Proposition \ref{p.roughness}).

\section{Proofs.}\label{sec:3}
\setcounter{section}{3}

Before proving Lemma \ref{l.isomorphism}, Lemma \ref{l.ED_diagonalized}, and Theorem \ref{thrm.ED_existence}, we begin by discussing, in detail, the spectral decomposition of the spherical Laplacian and the resulting multiplication symbols for the operators $A(r)$ and $A(r)^{1/2}$.

\subsection{Spectral Decomposition via Spherical Harmonics.}
\label{sec:spherical_harmonics}

For each degree 
$k\ge0$ with $k\in\Z$, there is an orthonormal set $\{Y_k^j\}_{j=1}^{\nu(k)}\subset L^2(S^{n-1},\C)$ which spans the $\nu(k)$-dimensional eigenspace of $\Delta_{S^{n-1}}$ and satisfies
$$
-\Delta_{S^{n-1}} Y_k^j = \mu_k Y_k^j \quad \text{where} \quad\mu_k= k(k+n-2).
$$

These $\{Y_k^j\}$ are the spherical harmonics of $\Delta_{S^{n-1}}$. For each 
$(u(r),v(r))\in  L^2(S^{n-1},\C^d)\times L^2(S^{n-1},\C^d)$
 they allow a spectral decomposition given by
\begin{align*}
	u(r,\varphi) & = \sum_{k = 0}^\infty\sum_{j=1}^{\nu(k)} u_{kj}(r)Y_k^j(\varphi), \quad u_{kj}(r) = \int_{S^{n-1}}u(r,\varphi)\overline{Y_k^j(\varphi)}d\varphi \in \C^d, \\
\quad	v(r,\varphi) & = \sum_{k = 0}^\infty\sum_{j=1}^{\nu(k)} v_{kj}(r)Y_k^j(\varphi), \quad v_{kj}(r) = \int_{S^{n-1}}v(r,\varphi)\overline{Y_k^j(\varphi)}d\varphi \in \C^d.
\end{align*}
Note that, unlike the Laplacian posed on $\R^n$, the spherical Laplacian posed on $S^{n-1}$ can be treated analogously in the case of $n=2$ and $n\ge 3$. If $n=2$, $\mu_k=k^2$ with $Y_k^1(\theta)=\exp({+\ri k\theta})$ and $Y_k^2(\theta)=\exp({-\ri k\theta})$ for $k>0$. See \cite{Frye2012sphericalharmonics} and the references therein.

This decomposition also allows for the definition of alternative, equivalent Sobolev norms on $H^\alpha(S^{n-1},\C^d)$, $\alpha>0$, and similarly on $L^2(S^{n-1},\C^d)$: 
\begin{equation*}\label{eq.decomp_norms}
	\|u(r,\cdot)\|_{H^\alpha}^2  := \sum_{k=0}^\infty\sum_{j=1}^{\nu(k)}(1+\mu_k)^\alpha\|u_{kj}(r)\|^2_{\C^d}, \qquad 
	\|u(r,\cdot)\|_{L^2}^2 := \sum_{k=0}^\infty\sum_{j=1}^{\nu(k)}\|u_{kj}(r)\|_{\C^d}^2.
\end{equation*}
According to this Fourier-like decomposition, \eqref{eq.principal_eq} 
can likewise be decomposed into $u_{kj}$ components  (suppressing the $j$ subscript, as each equation is independent of $j$)
\begin{equation*}\label{eq.spec_decomp_eq}
	\begin{pmatrix}
		\partial_r u_k \\
		\partial_r v_k
	\end{pmatrix} = \begin{pmatrix}
		0 & 1\\
		V_\infty + \mu_k/r^2 & 0
	\end{pmatrix}\begin{pmatrix}
		u_{k}\\v_{k}
	\end{pmatrix}
	=  \begin{pmatrix}
		0 & 1\\
		A_k(r) & 0
	\end{pmatrix}\begin{pmatrix}
		u_{k}\\v_{k}
	\end{pmatrix},
\end{equation*}
where $A_k(r) = V_\infty +\mu_k/r^2 \in \C^{d\times d}$ and $(u_k(r),v_k(r))\in\C^d\times \C^d$.

Notice that this spectral decomposition reveals $A(r)$ as a multiplication operator with symbol $A_k(r)$. Similarly, choosing the principal branch of the square root, we consider the operator $A(r)^{1/2}$ as a multiplication operator with symbol $A_k(r)^{1/2}$.

When $d>1$, we can further decompose the matrices $A_k(r)$  and $A_k(r)^{1/2}$ into scalar symbols via 
the diagonalization of $V_\infty$ given by hypothesis (H3). To do this, let $R$ be the matrix whose columns are given by the eigenvectors, $\{\rho_1, \cdots, \rho_d\}$, of $V_\infty$. 
 Then both $A_k(r)$ and $A_k(r)^{1/2}$ are also diagonalizable with the same eigenvectors as $V_\infty$:
\begin{align}
	A_k(r) 
	= R\begin{pmatrix}
		\lambda_1+\mu_k/r^2 & &\\[-1em]
		&\!\!\!\!\ddots &\\[-1em]
		& &  \!\!\!\!\lambda_d + \mu_k/r^2
	\end{pmatrix}R^{-1},
	\qquad
	A_k(r)^{1/2} 
	 = R \begin{pmatrix}
		\gamma_{k1}(r) & &\\[-1em]
		&\!\!\!\!\ddots &\\[-1em]
		& & \!\!\!\!\gamma_{kd}(r)
	\end{pmatrix}R^{-1}\notag
\end{align}
where $\gamma_{k\ell}(r) := \sqrt{\lambda_\ell + \mu_k/r^2}$, defined to be the principal square root, are the eigenvalues of $A_k(r)^{1/2}$.

In the proof of Lemma \ref{l.isomorphism}, we repeatedly use this Fourier-like decomposition
to obtain $r$-dependent upper-bounds for the $L^2(S^{n-1})$ norm of various operators related to $A(r)^{1/2}$. For example, for $A(r)^{1/2}$ we estimate
\begin{equation}   \label{eq.diag_norm}
	\begin{split}
	\Big\|A(r)^{1/2}u&\Big\|_{L^2(S^{n-1})}^2 \!
	=\sum_{k,j} \lp\|A_k(r)^{1/2}u_{kj}\rp\|^2_{\C^d}
	= \sum_{k,j} \lp\|R \begin{pmatrix}
		\gamma_{k1}(r) & &\\[-1em]
		&\!\!\!\!\ddots &\\[-1em]
		& & \!\!\!\!\gamma_{kd}(r)
	\end{pmatrix}R^{-1}u_{kj}\rp\|^2_{\C^d}\\
	\!&\le \sum_{k,j}\lp\|R\rp\|^2\lp(\max_{1\le \ell \le d}|\gamma_{k\ell}(r)|^2\rp)\lp\|R^{-1}\rp\|^2\lp\|u_{kj}\rp\|_{\C^d}^2 
		\le C\sum_{k,j} \max_{1\le\ell\le d}|\gamma_{k\ell}(r)|^2 \|u_{kj}\|_{\C^d}^2
	\end{split}
\end{equation}
where additional bounds on the scalar symbols $\gamma_{k\ell}(r)$ can then be applied. 

\begin{remark}
In the case of nontrivial Jordan blocks, a similar decomposition of the symbols $A_k(r)$ and $A_k(r)^{1/2}$ could be applied via the generalized eigenvectors of $V_\infty$. This decomposition would still reveal each $\gamma_{k\ell}(r)$ to be an eigenvalue of $A_k(r)^{1/2}$ corresponding to each distinct eigenvalue $\lambda_\ell$ of $V_\infty$.
For any $\lambda_\ell$ associated with a nontrivial Jordan chain, $\gamma_{k\ell}(r)$ would inherit the generalized eigenspace, though it would not inherit the same generalized eigenvectors. Consequently, for each nontrivial Jordan block of $R^{-1}V_\infty R$, the matrix $R^{-1} A_k(r)^{1/2} R$ would have a corresponding block of the same dimension which could be written as the sum of an $r$-dependent nilpotent matrix and the identity matrix multiplied by $\gamma_{k\ell}(r).$\\
Extending our results to this setting would require slightly different approaches to the proofs of Lemma \ref{l.isomorphism} and Lemma \ref{l.ED_diagonalized} in \S\ref{sec:lem_proofs}. For Lemma \ref{l.isomorphism}, employing a similar approach would require obtaining bounds on each nontrivial block analogous to the ones proven in the Appendix. Alternatively, the theory of primary matrix functions (see, Ch 6 of \cite{HornJohnson91}) suggests a less direct approach to proving the first statement of Lemma \ref{l.isomorphism}, but more investigation is required.  For Lemma \ref{l.ED_diagonalized}, though we give an explicit construction that relies on the diagonalizable nature of $V_\infty$, a less explicit proof still relying on robustness in $L^2\times L^2$ would likely follow even in the case of nontrivial Jordan blocks.
\end{remark}

\subsection{Proofs of Lemmas \ref{l.isomorphism} and \ref{l.ED_diagonalized}.}
\label{sec:lem_proofs}

In the proofs that follow we rely on a series of lemmas to bound scalar terms which are stated and proven in the Appendix. We also rely repeatedly on the spectral decomposition and diagonalization of $A(r)^{1/2}$ to decompose the $L^2(S^{n-1})$ norm as in \eqref{eq.diag_norm}.

\pagebreak[3]
\begin{proof}[Proof (of Lemma \ref{l.isomorphism}).] 
	
 We prove statement (i) in five steps:
	\begin{enumerate}[leftmargin=*,label=\arabic*.]
		\item \underline{$A(r)^{1/2}$ is well-defined and uniformly bounded in $r$.}\\[1em] Notice that there always exists a constant $K>0$, independent of both $\lambda_\ell$ and $r>1$, such that $|\gamma_{k\ell}(r)^2|=|\lambda_\ell+\mu_k/r^2|\le K(1+\mu_k)$. Then, by \eqref{eq.diag_norm}, we have
		\begin{align*}
			\|A(r)^{1/2}u\|_{L^2}^2 \le  \sum_{k,j} C \max_{1\le \ell\le d}|\gamma_{k\ell}(r)|^2\|u_{kj}\|^2_{\C^d} \le \sum_{k,j}\tilde{C}(1+\mu_k)\|u_{kj}\|^2_{\C^d} = \tilde{C}\|u\|_{H^1}^1. 
		\end{align*}
		\item \underline{$A(r)^{1/2}$ is strongly continuous.} \\[1em]
		Let $r_1,r_2\ge1$ be sufficiently close with $r_1\le 2r_2$. Then, again by our Fourier-like decomposition of the $L^2(S^{n-1})$ norm, as in \eqref{eq.diag_norm}, we have
			\begin{align*}
			\|A(r_1)^{1/2}u-A(r_2)^{1/2}u\|_{L^2}^2 &\le  \sum_{k,j}C\max_{1\le \ell \le d}|\gamma_{k\ell}(r_1)-\gamma_{k\ell}(r_2)|^2\|u_{k,j}\|^2_{\C^d}\\ 
			&\le \sum_{k,j}\tilde{C} |r_1-r_2|^2(1+\mu_k)\|u_{k,j}\|^2_{\C^d} \qquad\text{ by Lemma \ref{l:equiv_lem2}}.\\
			&= \tilde{C}|r_1-r_2|^2\|u\|_{H_1}^2 .
		\end{align*}
		\item \underline{$A(r)^{1/2}$ is differentiable.} \\[1em]
		We claim that the operator $\partial_r\!\lp[A(r)^{1/2}\rp]$, defined via the symbol $\partial_r\gamma_{k\ell}(r) = \left(\frac{-\mu_k/r^3}{\gamma_{k\ell}(r)}\right)$, is the Fr\'echet derivative of $A(r)^{1/2}$. To show this, fix $r\ge1$ and let $r'\ge1$ be sufficiently close to $r$ so that $r/2\le r'\le 2r$. Then, again by our Fourier-like decomposition, we have
			\begin{align*}
			\|A(r)^{1/2}u -& A(r')^{1/2}u - \partial_r A(r)^{1/2}(r-r')u\|_{L^2}^2\\ &\le \sum_{k,j}C\max_{1\le \ell\le d}|\gamma_{k\ell}(r) - \gamma_{k\ell}(r')-\partial_r\gamma_{k\ell}(r)(r-r')|^2\|u_{k,j}\|^2_{\C^d} \\
			&\le \tilde{C}|r-r'|^4\sum_{k,j}(1+\mu_k)\|u_{k,j}\|^2_{\C^d} \hspace{0.6in} \text{ by Lemma \ref{l:equiv_lem3}}.\\
			&= \tilde{C}|r-r'|^4\|u\|_{H^1}^2. 
		\end{align*}
		This leads us to conclude that $\partial_r A(r)^{1/2}$ satisfies the limit definition of the Fr\'echet derivative:
		\begin{align*}
			\lim_{r'\to r} \frac{\|A(r)^{1/2} - A(r')^{1/2}- \partial_r A(r)^{1/2}(r-r')\|_{L(L^2)}}{|r-r'|} \le \lim_{r'\to r}\frac{C'|r-r'|^2}{|r-r'|} = 0.
		\end{align*}
		\pagebreak[3]
		\item \underline{$\partial_rA(r)^{1/2}$ is continuous.}\\[1em]
		Let $r_1\ge r_2\ge 1$ be sufficiently close so that $r_2/2 \le r_1 \le 2r_2$. Then,  we have
		\begin{align*}
			\|\partial_rA(r_1)^{1/2}u-\partial_r&A(r_2)^{1/2}u\|_{L^2}^2 
			\le \sum_{k,j}C\max_{1\le \ell\le d}|\partial_{r}\gamma_{k\ell}(r_1)-\partial_r\gamma_{k\ell}(r_2)|^2\|u_{kj}\|^2_{\C^d} \\
			&\le \tilde{C}\left||r_1-r_2|+|r_1^3-r_2^3|\right|^2\sum_{k,j}(1+\mu_k)\|u_{k,j}\|^2_{\C^d} \quad \text{ by Lemma \ref{l:equiv_lem4}}\\
			&\le \tilde{C}\left(|r_1-r_2|+|r_1^3-r_2^3|\right)^2\|u\|_{H^1}^2.
		\end{align*}
		\item \underline{The operator norm of $\partial_r\!\left[A(r)^{1/2}\right]\!A(r)^{-1/2}: L^2\to L^2$ is of order $1/r$.} \\[1em]
		First note that the symbol of $A(r)^{-1/2}$ is given by $\gamma_{k\ell}(r)^{-1}$ and that, consequently, the symbol of the product $\partial_r\!\left[A(r)^{1/2}\right]\!A(r)^{-1/2}$ is given by $\left(\frac{-\mu_k/r^3}{\gamma_{k\ell}(r)^2}\right)$. We therefore obtain a bound on the $L^2$ operator norm of $\partial_r\!\left[A(r)^{1/2}\right]\!A(r)^{-1/2}$ by bounding the supremum of the symbols over all values of $k$ and $\ell$. 
		We find
		\begin{align*}
			\sup_{k,\ell} \left|\frac{-\mu_k/r^3}{\gamma_{k\ell}(r)^2}\right| = \left(\frac{1}{r}\right)\sup_{k,\ell}\frac{|\mu_k/r^2|}{|\lambda_\ell+\mu_k/r^2|} \le Cr^{-1},
		\end{align*}
		for some $C>0$ independent of $r$ by Hypothesis (H2).
		Consequently, by our Fourier-like decomposition of the $L^2(S^{n-1})$ norm, we have
		\begin{align*}
			\lp\|\partial_r\!\lp[A(r)^{1/2}\rp]\!A(r)^{-1/2}u\rp\|_{L^2}^2 & \le \sum_{k,j}C\max_{1\le \ell\le d} \lp|\frac{-\mu_k/r^3}{\gamma_{k\ell}(r)^2}\rp|\|u_{kj}(r)\|_{\C^d}^2 \\
			&\le \tilde{C}r^{-1}\sum_{kj}\|u_{kj}\|_{\C^d}^2 = \tilde{C}r^{-1}\|u\|_{L^2}^2,
		\end{align*}
		with $\tilde{C}$ also independent of $r$, so the operator norm on $\partial_r\!\lp[A(r)^{1/2}\rp]\!A(r)^{-1/2}$ is $O(r^{-1})$.
	\end{enumerate}
	
	\noindent To prove (ii), we first apply \eqref{eq.diag_norm} and estimate that
	\begin{align*}
		\|A(r)^{1/2}u\|_{L^2}^2 + \|v\|_{L^2}^2 &\le \sum_{k,j} C\max_{1\le \ell\le d} \left\{\left(|\lambda_\ell| + \mu_k/r^2\right)\right\}\|u_{k,j}\|^2_{\C^d} + \|v\|_{L^2}^2\\ 
		&\le \tilde{C}\left(\max_{1\le \ell\le d}\left\{|\lambda_\ell|\right\}\|u\|_{L^2}^2 + \frac{1}{r^2}\|u\|_{H^1}^2\right) + \|v\|_{L^2}^2 \le \tilde{C}\|\mbf{u}\|_{X_r}^2. 
	\end{align*}	
	On the other hand, 
	\begin{align*}
		\|\mbf{u}\|_{X_r}^2 &= \frac{1}{r^2}\|u\|_{H^1}^2 + \|u\|_{L^2}^2 + \|v\|_{L^2}^2\\
		&\le C\|u\|_{L^2}^2 + \|A(r)^{1/2}u\|_{L^2}^2 + \|v\|_{L^2}^2 \le \tilde{C}(\|A(r)^{1/2}u\|_{L^2}^2 +\|v\|_{L^2}^2)  \quad \text{ by Lemma \ref{l:equiv_lem5}}.
	\end{align*}
	This concludes the proof that $A(r)^{1/2}$ yields an isomorphism between $X_r$ and $L^2\times L^2$.
\end{proof}

The proof of Lemma \ref{l.ED_diagonalized} explicitly constructs an exponential dichotomy for \eqref{eq.diagonalized} from the hyperbolicity of $\mathcal{A}_\infty$ by exploiting the spectral decomposition and diagonalization of $A(r)^{1/2}$.

\begin{proof}[Proof (of Lemma \ref{l.ED_diagonalized}).] 
	
	Let $\tilde{\Phi}^{\rs/\ru}$ be defined as multiplication operators as follows:
		\begin{align}
		\tilde\Phi^\rs(r,t)\mbf{w} &= \begin{pmatrix}
			0 & 0 \\
			0 & \mathrm{exp}\lp(-\int_t^rA(s)^{1/2}ds\rp)
		\end{pmatrix}\mbf{w}\notag\\
		&=\sum_{k,j}  \begin{pmatrix}
			0 & 0 \\
			0 & \mathrm{exp}\lp(-\int_t^rA_k(s)^{1/2}ds\rp)
		\end{pmatrix}\mbf{w}_{kj}Y_{k}^j(\varphi),
		\quad r\ge t\notag\\
		\tilde\Phi^\ru(r,t)\mbf{w}  &= \begin{pmatrix}
			\mathrm{exp}\lp(\int_t^rA(s)^{1/2}ds\rp) & \,\,\,\,\,0 \\
			0 &  \,\,\,\,\,0
		\end{pmatrix}\mbf{w}\notag\\
		&=\sum_{k,j} \begin{pmatrix}
			\mathrm{exp}\lp(\int_t^rA_k(s)^{1/2}ds\rp) & \,\,\,\,\,0 \\
			0 &  \,\,\,\,\,0
		\end{pmatrix}\mbf{w}_{kj}Y_{k}^j(\varphi),
		\quad r\le t\label{eq.diag_ED}
	\end{align}
	where by the diagonalization of $V_\infty$, we can further decompose \eqref{eq.diag_ED} according to
	\begin{equation*}\label{eq.diag_ED_2}
		\mathrm{exp}\left(\pm\int_{t}^r A_k(s)^{1/2}ds \right)
		= R \begin{pmatrix}
			\mathrm{exp}\lp(\pm \int_t^r\gamma_{k1}(s)^{1/2}ds\rp)	& &\\[-1em]
			&\!\!\!\!\ddots &\\[-1em]
			& & \!\!\!\!\mathrm{exp}\lp(\pm\int_t^r\gamma_{kd}(s)^{1/2}ds\rp)	
		\end{pmatrix} R^{-1}.
	\end{equation*}
	The continuity and differentiability of \eqref{eq.diag_ED} follows from the continuity and differentiability of both $A(r)^{1/2}$ and each $\gamma_{k\ell}(r)$ (Lemmas \ref{l.isomorphism}, \ref{l:equiv_lem2}, \ref{l:equiv_lem3}, and \ref{l:equiv_lem4}).
	By taking the derivatives of $\tilde{\Phi}^{\rs}(r,t)\mbf{w}$ and $\tilde{\Phi}^{\ru}(r,t)\mbf{w}$ with respect to $r$ one can directly confirm that for each $t\in\R$ and $\mbf{w}\in (L^2(S^{n-1})\times L^2(S^{n-1}))$ these are both solutions to \eqref{eq.diagonalized}.  
	The invariance conditions required by exponential dichotomies (as in Definition \ref{d.ed-Xr})
	are also clearly satisfied by the evolutionary operators and their corresponding projections,  $\tilde{P}^\rs:=\tilde{\Phi}^\rs(r,r) = \lp(\begin{smallmatrix}0 & 0\\ 0& 1\end{smallmatrix}\rp)$ and $\tilde{P}^\ru:=\tilde{\Phi}^\ru(r,r) =  \lp(\begin{smallmatrix}1 & 0\\ 0& 0\end{smallmatrix}\rp)$.
	
	It remains to be shown that the exponential decay estimates hold for 
	$$\tilde{\eta} := \min_{1\le\ell\le d}\mathrm{Re}[\lambda_\ell^{1/2}].$$
	These decay estimates follow from the fact that for all $k\in\{0,1,2,\ldots\}$ and $r>0$, the following bound holds: $\mathrm{Re}\lp[\gamma_{k\ell}(r)\rp]\ge \mathrm{Re}[\lambda_{\ell}^{1/2}]$ (see Lemma \ref{l.real_lower_bd}). We therefore obtain:
	\begingroup\allowdisplaybreaks[3]
		\begin{align*}
			\left\|\tilde{\Phi}^\rs(r,t)\mbf{w}\right\|_{L^2\times L^2}^2 
			&= \sum_{k,j}\left\|\mathrm{exp}\left(-\int_t^rA_k(s)^{1/2}ds\right)
			w^-_{kj}\right\|_{\C^d}^2\\
			&= 
			\sum_{k,j}\left\|R \begin{pmatrix}
				\mathrm{exp}\lp(- \int_t^r\gamma_{k1}(s)^{1/2}ds\rp)	& &\\[-1em]
				&\!\!\!\!\ddots &\\[-1em]
				& & \!\!\!\!\mathrm{exp}\lp(-\int_t^r\gamma_{kd}(s)^{1/2}ds\rp)	
			\end{pmatrix} R^{-1}w^{-}_{kj}\right\|_{\C^d}^2\\
			&\le\sum_{k,j}\|R\|^2
			\left(\max_{1\le \ell \le d}\left| \mathrm{exp}\left(-\int_t^r\gamma_{k\ell}(s)ds\right) \right|^2\right)
			\|R^{-1}\|^2\|w^-_{kj}\|^2_{\C^d}\\
			&\le \sum_{k,j}\|R\|^2\|R^{-1}\|^2
			\left(\max_{1\le \ell \le d}
			\left[\mathrm{exp}\left({-\mathrm{Re}\big[\lambda_\ell^{1/2}\big](r-t)}\right)\right]^2
			\right)\|w^-_{kj}\|_{\C^d}^2\\
			&\le \sum_{k,j}C
			\left(\mathrm{exp}\left({-\min\limits_{1\le \ell \le d}2\mathrm{Re}\big[\lambda_\ell^{1/2}\big](r-t)}\right)
			\right)\|w^-_{kj}\|_{\C^d}^2\\
			&\le \sum_{k,j}C \exp({-2\tilde{\eta}(r-t)})\|w^-_{kj}\|_{\C^d}^2
			= C\exp({-2\tilde{\eta}(r-t)})\|w^-\|_{L^2}^2\\
			&\le C\exp({-2\tilde{\eta}(r-t)})\|\mbf{w}\|_{L^2\times L^2}^2, \qquad r\ge t
		\end{align*}
		\endgroup
		A nearly identical calculation shows the analogous statement for $\tilde{\Phi}^\ru(r,t)$ on $r\le t$.
	\end{proof}

\subsection{Proof of Theorem \ref{thrm.ED_existence}.}
\label{sec:thrm_proof}

\begin{proof}[Proof (of Theorem \ref{thrm.ED_existence}).] 
	By Lemma \ref{l.ED_diagonalized}, \eqref{eq.diagonalized} has an exponential dichotomy in $L^2(S^{n-1})\times L^2(S^{n-1})$ for $r>1$ with decay rate $\tilde{\eta}$. Recall that \eqref{eq.diagonalized} is related to the principal part of \eqref{eq.COV_eq} via a change of coordinates $w^\pm = \tilde{u}\pm v$. This change of coordinates is equivalent to conjugation by the linear transformation $M:= \left(\begin{smallmatrix}
		1 & \,\,\,1\\ 1& -1
	\end{smallmatrix}\right)$:
	\[\mbf{w}_r = M \begin{pmatrix}
		0& A(r)^{1/2}\\  A(r)^{1/2} &0
	\end{pmatrix}\! M^{-1} \,\mbf{w} \qquad \text{and}\qquad \mbf{\tilde u}_r = M^{-1} \begin{pmatrix} A(r)^{1/2} &0  \\0 & -A(r)^{1/2}\end{pmatrix}\! M\,\mbf{\tilde u}.\]
	Notice that since $M = 2M^{-1}$, the related transformation given by $\tilde{M} := \tfrac{1}{\sqrt{2}} M$ is a unitary, self-adjoint operator. Consequently, conjugation by $\tilde{M}$ (equivalent to conjugation by $M$) defines an  isometric isomorphism on $L^2 \times L^2$, and the principal part of \eqref{eq.COV_eq} has an exponential dichotomy on $L^2 \times L^2$ for $r>1$ with decay rate $\tilde{\eta}$ defined by
	\begin{align}
		\hat{\Phi}^\rs(r,t) &:= M^{-1}\tilde{\Phi}^{\rs}(r,t)M, \qquad r\ge t\notag\\
		\hat{\Phi}^\ru(r,t) &:= M^{-1}\tilde{\Phi}^{\ru}(r,t)M, \qquad r\le t\label{eq.principal_COV_ED}
	\end{align}
	where $\tilde{\Phi}^{\rs/\ru}$ give the exponential dichotomy for \eqref{eq.diagonalized} defined by \eqref{eq.diag_ED}.
	
	From \eqref{eq.principal_COV_ED} we can directly apply the robustness of exponential dichotomies (Theorem 1 in \cite{PeterhofSandstedeScheel97})  to get an exponential dichotomy $\hat\Phi^{\rs/\ru}_{*}$ for \eqref{eq.COV_eq} in $L^2 \times L^2$ for $r>1$ again with the same decay rate $\tilde \eta$, noting that the perturbed part of \eqref{eq.COV_eq} is $\mathcal{O}(1/r)$. 
	
	Next, we use the exponential dichotomy $\hat\Phi^{\rs/\ru}_*$ for \eqref{eq.COV_eq} on $L^2\!\!\times\! L^2$ and the isomorphism established by Lemma \ref{l.isomorphism} to construct an exponential dichotomy for \eqref{eq.principal_eq} on $X_r$.
	
	Notice that by Lemma \ref{l.isomorphism} \textit{(ii)} the transformation $(u,v)\mapsto (A(r)^{1/2}u,v)$ and can be written as 
	\[T(r) := \begin{pmatrix}A(r)^{1/2} & 0\\ 0 & 1\end{pmatrix}: X_r\to L^2\!\!\times\! L^2 \quad \text{with} \quad T(r)^{-1} =  \begin{pmatrix}A(r)^{-1/2} & 0\\ 0 & 1\end{pmatrix}: L^2\!\!\times\! L^2\to X_r.\]
	Therefore, \eqref{eq.COV_eq} is equivalent to $\mbf{\tilde u}_r = \big[T(r)\tilde{\mathcal{A}}(r) T(r)^{-1} + \partial_r\lp[T(r)\rp]T(r)^{-1}\big]\mbf{\tilde u}$. 
	
	A similar $r$-dependent conjugation by $T$ will transform  $\hat\Phi^{\rs}_*(r,t)$ and $\hat\Phi^\ru_*(r,t)$ from a family of solution operators for \eqref{eq.COV_eq} in $L^2\times L^2$ to a family of solution operators for \eqref{eq.principal_eq} in $X^r$: 
	\begin{align}
		\Phi^{\rs}_*(r,t) &:= T(r)^{-1} \hat\Phi_*^{\rs}(r,t)T(t): X_t \to X_r, \qquad r\ge t\notag\\
		\Phi^{\ru}_*(r,t) &:= T(r)^{-1} \hat\Phi_*^{\ru}(r,t)T(t): X_t \to X_r, \qquad r\le t.  \label{eq.principal_ED}
	\end{align}
	
	The continuity and differentiability of $\Phi^{\rs/\ru}_*$ follow directly from Lemma \ref{l.isomorphism} and the continuity and differentiability of $\hat{\Phi}^{\rs/\ru}_*$. The projection and invariance properties required by Definition \ref{d.ed-Xr} are also inherited from $\hat{\Phi}^{\rs/\ru}_*$. All that remains to show is that $\Phi^{\rs}_*$ and $\Phi^{\ru}_*$ satisfy exponential decay with rate $\tilde{\eta}$; this is addressed by the second statement of Lemma \ref{l.isomorphism}.
	
	 Though the isomorphism $T(r)$ is not an isometry, by Lemma \ref{l.isomorphism} it does respect the norm up to constants independent of $r$. Therefore the family $\Phi_*^{\rs/\ru}$ inherits the same exponential decay properties in $X_r$ as the exponential dichotomy $\hat\Phi_*^{\rs/\ru}$ has in $L^2\!\!\times\! L^2$, including the same decay rate $\tilde\eta$. For any $\mbf{u}\in X_t$,
	 \begin{align*}
	 	\lp\|\Phi_*^\rs(r,t) \mbf{u}\rp\|_{X_r} 
	 	\le C_2\lp\|\hat\Phi_*^\rs(r,t)T(t)\mbf{u}\rp\|_{L^2\!\times L^2}
	 	&\le C_2K\exp({-\tilde\eta(r-t)})\lp\|T(t)\mbf{ u}\rp\|_{L^2\!\times L^2}\\
	 	& \le \frac{C_2K}{C_1}\exp({-\tilde\eta(r-t)})\lp\|\mbf{u}\rp\|_{X_t}, && r\ge t,\\[1em]
	 	\lp\|\Phi_*^\ru(r,t) \mbf{u}\rp\|_{X_r}
	 	\le C_2\lp\|\hat\Phi_*^\ru(r,t)T(t)\mbf{ u}\rp\|_{L^2\!\times L^2}
	 	&\le C_2K\exp({\tilde\eta(r-t)})\lp\|T(t)\mbf{ u}\rp\|_{L^2\!\times L^2}\\
	 	& \le \frac{C_2K}{C_1}\exp({\tilde\eta(r-t)})\lp\|\mbf{u}\rp\|_{X_t}, && r\le t.
	 \end{align*}
	 
	Therefore \eqref{eq.principal_ED} defines an exponential dichotomy for the principal equation \eqref{eq.principal_eq} in $X_r$ for $r>1$ with decay rate $\tilde\eta$. Proposition \ref{p.roughness} applied to \eqref{eq.principal_eq} with perturbation $\tilde{B}(r)\in \mathcal{O}(1/r)$ from \eqref{eq.perturb_op} gives the final result: $\Phi_*^{\rs/\ru}$ can be perturbed to an exponential dichotomy $\Phi^{\rs/\ru}$ for \eqref{eq.spatial_eq} with decay rate $\eta = \tilde\eta$. 	
\end{proof}

\section{Conclusion.}\label{sec:4}

Theorem \ref{thrm.ED_existence}, together with the results of \cite{Beck21-ExpDich}, aid in the characterization of the bounded behavior of spatial evolutionary systems that arise from elliptic PDEs of the form  $\Delta u+ V(x) u=0$ where $V(x)=V(r,\cdot)$ decays to a constant potential sufficiently fast as $r$ goes to infinity. In fact, exponential dichotomies provide a construction of bounded solutions to nonlinearly perturbed equations via a Duhamel formula. This gives rise to a constructive approach to proving the existence of localized stationary solutions of pattern-forming PDEs that lack any kind of constraining symmetry. Radially symmetric localized solutions were thoroughly characterized by \cite{Scheel03}. Recent work in proving the existence of localized patterns with dihedral symmetries (e.g. hexagonal patterns) utilized a combination of analytic and numerical techniques; see \cite{BramburgerHillLloyd25} for a comprehensive review. In \cite{Cadiot25}, computer-assisted proof techniques via Newton-Kantorovich methods have recently proven successful for rigorously verifying the existence of localized stationary solutions to Swift-Hohenburg without restricting to specific symmetries. However, with our contribution of far-field exponential dichotomies in the setting of unbounded multidimensional spatial domains, we hope to open up the tools of spatial dynamics as set forward by \cite{PeterhofSandstedeScheel97} and \cite{SandstedeScheel04_classify} to this new setting.

To see the broad applicability of Theorem \ref{thrm.ED_existence}, notice that any stationary reaction diffusion model of the form $0=\Delta u + f(u)$, $u\in\C^d$, becomes $0=\Delta u + Df(q(x))u$ when linearized about a stationary solution $q(x)$. Notice also that this matrix  $Df(q(x))$ will certainly decay to the constant coefficient matrix $Df(0)$ as $|x|\to\infty$ if the solution $q(x)$ is localized. Therefore the linearization of such a reaction diffusion stationary problem on $\R^n$ about any localized solution $q(x)$ has exponential dichotomies in both the core field and the far field provided that $Df(0)$ has eigenvalues off of the negative real line. 

A natural extension would be to perturb the reaction diffusion system with a small inhomogenous nonlinearity,
\begin{equation}\label{eq.perturbed_RD}
	0=\Delta u + f(u) + \mu h(x,u,\mu) = \Delta u + Df(q(x))u +
	\big[f(u) - Df(q(x))u +\mu h(x,u,\mu) \big],
\end{equation}
and use the core-field and far-field dichotomies for the linearization to construct bounded solutions in both the core field ($r\in(0,1]$) and the far field ($r\in [1,\infty)$). Further analysis from the channel-setting, for example, may be extended to match these solutions and prove the existence of a bounded solution on the whole domain.
Of interest would be to characterize how general the inhomogeneity $h$ can be while preserving smooth perturbation of $q(x)$ and absence of essential spectrum.

A similar construction was demonstrated in the channel setting, $(x,y)\in\R\times\Omega$ with $\Omega$ compact in $\R^{n-1}$, by \cite{PeterhofSandstedeScheel97} for a spatial evolutionary equation of the form:
\begin{equation}\label{eq.PSS_perturb}
	u_x = Au + G(u)+\mu H(x,u,\mu), \qquad u\in X=L^2(\Omega), 
\end{equation}
where $A$ is an autonomous, linear (differential) operator on $L^2(\Omega)$ and $H$ is periodic in $x$. Assuming that a homoclinic solution $q(x)$ exists for $\mu=0$ in \eqref{eq.PSS_perturb} and that the linearization about $q$ has an exponential dichotomy on $\R^+$ and on $\R^-$, \cite{PeterhofSandstedeScheel97} constructed a bounded solution on $\R^+$ and $\R^-$, respectively. Then, by applying a Melnikov argument with Lyapunov-Schmidt techniques, \cite{PeterhofSandstedeScheel97} showed that the stable space on $\R^+$ and the unstable space on $\R^-$ intersect in such a way that a unique bounded solution exists on the whole domain $\R$ when $\mu>0$ is small.

Ongoing and future work aims to extend such an argument to our more general domain with the primary challenge being the higher dimensional kernel that comes from having not just one translational symmetry but  at least $n$ symmetries, one for each unbounded direction in state space. With such an argument, the existence of exponential dichotomies in the setting described above would allow us to take previously understood radially symmetric localized solutions and perturb them in non-radially symmetric ways to prove the existence of arbitrary localized stationary patterns. 

\needspace{3\baselineskip}
It would also be of interest to investigate the connection of our approach to the general works on solutions of singularly perturbed nonlinear elliptic problems of the form $-\epsilon^2\Delta u + V(|x|)u - u|u|^{p-1} = 0,$ posed on $\R^n$ for suitable $p$, where radial solutions concentrated on spherical subsets related to minima of a certain modified potential have been identified, along with an sequence of non-radial solutions bifurcating as $\epsilon\rightarrow0^+$; see \cite{ambrosetti_singularly_2003, ambrosetti2006perturbation} and references therein.

\appendix
\section*{Appendix}
\label{sec:Apdx}

\renewcommand{\theLemmaApdx}{A.\arabic{LemmaApdx}}

\begin{LemmaApdx}\label{l:equiv_lem0} For all $k\ge0$, $\mu_k^{1/2} \le (1+\mu_k)^{1/2}$. Likewise, for $r\ge 1$, $\mu_k^{1/2}/r\le (1+\mu_k)^{1/2}$. Furthermore, $$\displaystyle\sum_{k,j}\mu_k\|u_{kj}\|^2_{\C^d} \le \|u\|_{H^1}^2.$$ 
\end{LemmaApdx}

\begin{remark}
	Though Lemma \ref{l:equiv_lem0} is trivial, since $\mu_k\ge 0$ is real valued, we find it convenient to refer to this bound within long calculations for the ease of the reader.
\end{remark}

\begin{LemmaApdx}\label{l:equiv_lem1}
	Given \textup{(H2)}, for each $k\ge0$, $1\le\ell\le d$, and any arbitrary $r_1,r_2 \ge 1$,
	\[|\gamma_{k\ell}(r_1)+\gamma_{k\ell}(r_2)|\ge |\gamma_{k\ell}(r_j)|,\qquad j=1,2.\]
\end{LemmaApdx}

\begin{proof}
	Recall that $\gamma_{k\ell}(r)^2 = \lambda_\ell + \mu_k/r^2$ where $\mu_k\ge 0$ and $\lambda_\ell\in\C\backslash(-\infty,0]$. If $\Im\lp[\lambda_\ell\rp] \ge 0$, then by taking the principal branch of the square root, we also have $\Im\lp[\gamma_{k\ell}(r)\rp]\ge 0$ for all $r\ge 1$. Similarly, if $\Im\lp[\lambda_\ell\rp]\le 0$, then $\Im\lp[\gamma_{k\ell}(r)\rp]\le 0$ for all $r\ge 1$. Since for each $\ell$ the imaginary parts of $\gamma_{k\ell}(r_1)$ and $\gamma_{k\ell}(r_2)$ have matching sign, we have that
	\[\Im\left[\gamma_{k\ell}(r_1)+\gamma_{k\ell}(r_2)\right]^2 = \left(\Im\left[\gamma_{k\ell}(r_1)\right]+\Im\left[\gamma_{k\ell}(r_2)\right]\right)^2 \ge \Im\left[\gamma_{k\ell}(r_j)\right]^2.\]
	Furthermore, because we are taking the principal branch of the square root, $\Re\left[\gamma_{k\ell}(r)\right] > 0$ for all $r$ and all $\ell$. Therefore we have
		\[\Re\left[\gamma_{k\ell}(r_1)+\gamma_{k\ell}(r_2)\right]^2 = \left(\Re\left[\gamma_{k\ell}(r_1)\right]+\Im\left[\gamma_{k\ell}(r_2)\right]\right)^2 \ge \Re\left[\gamma_{k\ell}(r_j)\right]^2.\]
	Consequently, we have that
	\begin{align*}
		\lp|\gamma_{k\ell}(r_1)+\gamma_{k\ell}(r_2)\rp| 
			&= \sqrt{\Re\left[\gamma_{k\ell}(r_1)+\gamma_{k\ell}(r_2)\right]^2+\Im\left[\gamma_{k\ell}(r_1)+\gamma_{k\ell}(r_2)\right]^2 }\\
			&\ge \sqrt{\Re\left[\gamma_{k\ell}(r_j)\right]^2+\Im\left[\gamma_{k\ell}(r_j)\right]^2} 
			=  |\gamma_{k\ell}(r_j)|.
	\end{align*}
\end{proof}

\begin{LemmaApdx}\label{l:equiv_lem1.5}
	Given \textup{(H2)}, for $k>0$, $1\le\ell\le d$, and any $r_1,r_2 \ge 1$,
	\[\left|\frac{\mu_k^{1/2}/r_j}{\gamma_{k\ell}(r_1)+\gamma_{k\ell}(r_2)}\right| \le C, \qquad j=1,2\]
	where $C>0$ is some constant independent of $k$, $\ell$, and $r_{1,2}$.
\end{LemmaApdx}

\begin{proof}
	By Lemma \ref{l:equiv_lem1}, we have that
	\begin{align*}
		\left|\frac{\mu_k^{1/2}/r_j}{\gamma_{k\ell}(r_1)+\gamma_{k\ell}(r_2)}\right|^2 
				\le \frac{\mu_k/r_j^2}{|\gamma_{k\ell}(r_j)|^2} = \frac{\mu_k/r_j^2}{\big|\gamma_{k\ell}(r_j)^2\big|}
				 = \frac{\mu_k/r_j^2}{\big|\lambda_\ell+\mu_{k}/r_j^2\big|} = f_\ell\left(\mu_k/r_j^2\right).
	\end{align*}
	where $f_\ell(\xi):= \dfrac{\xi}{|\lambda_\ell+\xi|}$. 
	
	Notice that $f_\ell(\xi)$ is well-defined and continuous on the domain $\xi\in [1,\infty)$ since $\lambda_\ell \in\C\backslash(-\infty,0]$, and that $\lim_{\xi\to\infty}f(\xi) = 1$.  Furthermore, it is straightforward to show using calculus techniques that $f_\ell(\xi)$ is strictly increasing for $\xi\ge1$ provided $\Re\lp[\lambda_\ell\rp]\ge 0$. If $\Re\lp[\lambda_\ell\rp]<0$, then by solving $f'_\ell(\xi)=0$ we find that $f_\ell(\xi)$  attains a single maximum on $[1,\infty)$ at $\xi=-|\lambda_\ell|^2/\Re\lp[\lambda_\ell\rp]$ given by
	\[f_\ell\left(\frac{|\lambda_\ell|^2}{\lp|\Re\lp[\lambda_\ell\rp]\rp|}\right) = \frac{|\lambda_\ell|}{\lp|\Im\left[\lambda_\ell\right]\rp|}.\]
	
	The result follows with constant $C>0$ chosen as
	\[C=\max_{1\le \ell\le d}\left(\sup_{\xi\ge 1} \sqrt{f_\ell(\xi)}\right).\]
\end{proof}

\begin{LemmaApdx}\label{l:equiv_lem2}
	Given \textup{(H2)}, for each $k>0$, $1\le \ell\le d$, and $r_1\ge r_2\ge 1$ sufficiently close so that $r_1\le 2r_2$,
	$$|\gamma_{k\ell}(r_1) - \gamma_{k\ell}(r_2)| \le C(1+\mu_k)^{1/2}|r_1-r_2|.$$		
\end{LemmaApdx}

\begin{proof}
	We begin by multiplying $\gamma_{k\ell}(r_1)-\gamma_{k\ell}(r_2)$ by $(\gamma_{k\ell}(r_1)+\gamma_{k\ell}(r_2))$ in both the numerator and the denominator, and simplifying the result into a single fraction:
	\begin{align*}
		| \gamma_{k\ell}(r_1) - \gamma_{k\ell}(r_2)|&= \lp| \frac{\mu_k(r_1+r_2)(r_1-r_2)}{r_1^2r_2^2(\gamma_{k\ell}(r_1)+\gamma_{k\ell}(r_2))}\rp|\notag\\
		&\leq |r_1 - r_2|\left|\frac{r_1+r_2}{r_1^2r_2^2}\right|\lp|\frac{\mu_k}{(\gamma_{k\ell}(r_1) + \gamma_{k\ell}(r_2))}\rp|.\notag
	\end{align*}
	We next apply 
	 our assumption that $r_2\le r_1\le 2r_2$ to bound the coefficients via the estimates $r_1+r_2\le 3 r_2$ and $r_1^2r_2^2\ge r_2^4$. For $C=3$, this yields
	\begin{align*}
		| \gamma_{k\ell}(r_1) - \gamma_{k\ell}(r_2)|
		&\le C|r_1 - r_2|\lp|
										\frac{\mu_k/r_2^3}
												{\gamma_{k\ell}(r_1)+\gamma_{k\ell}(r_2)}\rp|.
	\end{align*}
	By factoring $\mu_k/r_2^3$ into $(\mu_k^{1/2}/r_2^2)\cdot(\mu_k^{1/2}/r_2)$ and applying Lemmas \ref{l:equiv_lem0} and \ref{l:equiv_lem1.5}, we then obtain
	\begin{align*}
		| \gamma_{k\ell}(r_1) - \gamma_{k\ell}(r_2)|&\le C|r_1 - r_2|\lp|\mu_k^{1/2}/r_2^2\rp| \le C |r_1-r_2| (1+\mu_k)^{1/2}.
	\end{align*}
\end{proof}

\begin{LemmaApdx}\label{l:equiv_lem3}
	Let $r_1,r_2>0$ be sufficiently close so that $r_1/2\le r_2 < 2r_1$. Then, given \textup{(H2)}, for all $k>0$ and $1\le \ell\le d$,
	$$|\gamma_{k\ell}(r_1) - \gamma_{k\ell}(r_2) - \partial_r\gamma_{k\ell}(r_1)(r_1-r_2)| \le C(1+\mu_k)^{1/2}|r_1-r_2|^2.$$		
\end{LemmaApdx}

\begin{proof}
	The result follows from a rather tedious calculation in which we suppress the $\ell$ subscript. Notice that in the first step of the calculation we simplify $(\gamma_{k\ell}(r_1)-\gamma_{k\ell}(r_2))$ as in the proof of Lemma \ref{l:equiv_lem2}.
	\begin{align*}
		\gamma_k(r_1) - &\gamma_k(r_2) - \p_r\gamma_k(r_1)(r_1 - r_2)
		= \frac{\mu_k(r_2+r_1)(r_2-r_1)}{r_1^2r_2^2(\gamma_{k}(r_1)+\gamma_{k}(r_2))}
		- \frac{\mu_k(r_2 - r_1)}{r_1^3 \gamma_k(r_1) }\notag\\[4pt]
		&= (r_2 - r_1)\frac{\mu_k}{r_1^2}\lp[\frac{r_2+r_1}{r_2^2(\gamma_k(r_1)  + \gamma_k(r_2) )} - \frac{1}{r_1\gamma_k(r_1) } \rp]\notag\\[4pt]
		&= (r_2 - r_1)\frac{\mu_k}{r_1^2}\lp[\frac{r_1(r_1+r_2)\gamma_k(r_1)  - r_2^2(\gamma_k(r_1)  + \gamma_k(r_2) )}{r_1r_2^2\gamma_k(r_1) (\gamma_k(r_1)  + \gamma_k(r_2)) }\rp]\notag\\[4pt]
		&= (r_2 - r_1)\frac{\mu_k}{r_1^2}\lp[\frac{ r_1^2(\gamma_k(r_1)  - \gamma_k(r_2) ) + \gamma_k(r_2) (r_1^2-r_2^2)+r_2\gamma_k(r_1) (r_1 - r_2) }{r_1r_2^2\gamma_k(r_1) (\gamma_k(r_1)  + \gamma_k(r_2) )}\rp]   
	\end{align*}
	We pause to note that this last equality is obtained by adding and subtracting the term $r^2_1\gamma_{k}(r_2)$. Continuing with the calculation, 
	this last term is equivalent to:
	\begin{align*}
		(r_2-r_1)&\frac{\mu_k}{r_1^2}\Bigg[\frac{\mu_k(r_2+r_1)(r_2-r_1)}{r_1r_2^4\gamma_k(r_1)(\gamma_{k}(r_1)+\gamma_{k}(r_2))^2}
		+(r_1-r_2)\lp(\frac{\gamma_k(r_2)(r_1+r_2)+r_2\gamma_k(r_1)}{r_1r_2^2\gamma_k(r_1)(\gamma_k(r_1)+\gamma_k(r_2))}\rp)
		\Bigg]\notag\\
		=& (r_2 - r_1)^2\frac{\mu_k}{r_1^2}\Bigg[ \frac{\mu_k(r_2+r_1)}{r_1r_2^4 \gamma_k(r_1) (\gamma_k(r_1)  + \gamma_k(r_2) )^2} \\
		&\hspace{1.1in}
		- \frac{\gamma_k(r_2) (r_1+r_2)}{r_1r_2^2\gamma_k(r_1) (\gamma_k(r_1)  + \gamma_k(r_2) )} 
		-\frac{1}{r_1r_2(\gamma_k(r)  + \gamma_k(r_2) )} \Bigg]\notag\\
		=& (r_2-r_1)^2\frac{\mu_k^{1/2}}{r_1} \Bigg[\frac{r_1(r_2+r_1)}{r_2^4}\frac{\mu_k^{3/2}/r_1^3}{\gamma_k(r_1) (\gamma_k(r_1)  + \gamma_k(r_2) )^2} \\ 
		&\hspace{1.1in}-\frac{(r_2+r_1)}{r_1r_2^2}\frac{\gamma_k(r_2)  (\mu_k^{1/2}/r_1)}{\gamma_k(r_1) (\gamma_k(r_1)  + \gamma_k(r_2) )}
		-\frac{1}{r_1 r_2}\frac{\mu_k^{1/2}/r_1}{(\gamma_k(r_1)  + \gamma_k(r_2))} \Bigg]
	\end{align*} 
	From this calculation we have that
	\begin{align*}
		\big|\gamma_k(r_1) - &\gamma_k(r_2) - \p_r\gamma_k(r_1)(r_1 - r_2)\big| \le (r_1-r_2)^2\frac{\mu_k^{1/2}}{r_1}\Big[|I|+|II|+|III|\Big]
	\end{align*}
	where, by Lemma \ref{l:equiv_lem1.5}, the assumption that $r_1/2\le r_2\le 2r_1$, and the continuity of $\gamma_k$ on $r_1\ge 1$, we have
		\begin{align*}
		|I| =
				 \left(\frac{r_1(r_2+r_1)}{r_2^4}\right) \lp|\frac{\left(\mu_k^{1/2}/r_1\right)^3}{\gamma_k(r_1) (\gamma_k(r_1)  + \gamma_k(r_2) )^2}\rp|
			&\le 3\cdot 2^4 \,\lp|\frac{\left(\mu_k^{1/2}/r_1\right)^3}{\gamma_k(r_1) (\gamma_k(r_1)  + \gamma_k(r_2) )^2}\rp|\leq C\\
		|II| = 
				\left(\frac{r_2+r_1}{r_1r_2^2}\right)\lp|\frac{\gamma_k(r_2) }{\gamma_k(r_1) }\rp| \lp| \frac{\mu_k^{1/2}/r_1}{(\gamma_k(r_1)  + \gamma_k(r_2) )}\rp|
			&\leq 3\cdot4 \lp|\frac{\gamma_k(r_2) }{\gamma_k(r_1) }\rp| \lp| \frac{\mu_k^{1/2}/r_1}{(\gamma_k(r_1)  + \gamma_k(r_2) )}\rp|\leq C\\
		|III| =\quad
				\left(\frac{1}{r_1r_2}\right) \quad \lp| \frac{\mu_k^{1/2}/r_1}{(\gamma_k(r_1)  + \gamma_k(r_2) )}\rp|
		\quad &\le \,\,2\,\,\lp| \frac{\mu_k^{1/2}/r_1}{(\gamma_k(r_1)  + \gamma_k(r_2) )}\rp|\leq C.
	\end{align*}
	Finally, the simple bound in Lemma \ref{l:equiv_lem0}, namely, $\mu_k^{1/2}/r_1 \le (1+\mu_k)^{1/2}$, concludes the proof.
\end{proof}

\begin{LemmaApdx}\label{l:equiv_lem4}
	Given \textup{(H2)}, for any arbitrary $r_1\ge r_2 \ge 1$ with $r_1$ and $r_2$ sufficiently close so that $r_2/2\le r_1\le 2r_2$, 
	$$|\partial_r\gamma_{k\ell}(r_1) - \partial_r\gamma_{k\ell}(r_2)| \le C(1+\mu_k)^{1/2}\left(|r_1-r_2|+|r_1^3-r_2^3|\right).$$		
\end{LemmaApdx}

\begin{proof}
	We begin by writing the difference as a single fraction. 
	\begin{align*}
		|\partial_r \gamma_{k\ell}(r_1) - \partial_r \gamma_{k\ell}(r_2)| &= \mu_k \left|\frac{1}{r_2^3\gamma_{k\ell}(r_2)}-\frac{1}{r_1^3\gamma_{k\ell}(r_1)}\right|
		=\mu_k \left|\frac{r_1^3\gamma_{k\ell}(r_1)-r_2^3 \gamma_{k\ell}(r_2)}{r_1^3 r_2^3\gamma_{k\ell}(r_2)\gamma_{k\ell}(r_1)}\right|\\[4pt]
		&=\mu_k\left|\frac{r_1^3(\gamma_{k\ell}(r_1)-\gamma_{k\ell}(r_2))-\gamma_{k\ell}(r_2)(r_2^3-r_1^3)}{r_1^3 r_2^3 \gamma_{k\ell}(r_2)\gamma_{k\ell}(r_1)}\right|
	\end{align*}
	The second line of the calculation above is obtained by adding and subtracting $r_1^3\gamma_{k\ell}(r_2)$ in the numerator. Breaking up this fraction and applying triangle inequality we obtain:
	\begin{align*}
		|\partial_r \gamma_{k\ell}(r_1) - \partial_r \gamma_{k\ell}(r_2)| &\le \left(\frac{\mu_k}{r_2^3|\gamma_{k\ell}(r_2)\gamma_{k\ell}(r_1)|}\right)\!\left|\gamma_{k\ell}(r_1)-\gamma_{k\ell}(r_2)\right| 
		+\! \left(\frac{\mu_k}{r_1^3r_2^3|\gamma_{k\ell}(r_1)|}\right) \left|r_2^3 -r_1^3\right|.
	\end{align*}
	Applying the same argument as in Lemma \ref{l:equiv_lem1.5} we obtain
	\begin{align*}
			\frac{\mu_k^{1/2}/r_j}{|\gamma_{k\ell}(r_j)|}\le C, \qquad j=1,2.
	\end{align*} 
	Furthermore, by our assumption that $r_2/2\le r_1 \le 2r_2$, we have 
	\[\frac{\mu_k/r_2^2}{|\gamma_{k\ell}(r_1)\gamma_{k\ell}(r_2)|} \le \left(\frac{2\mu_k^{1/2}/r_1}{|\gamma_{k\ell}(r_1)|}\right)\left(\frac{\mu_k^{1/2}/r_2}{|\gamma_{k\ell}(r_2)|}\right) \le C.\]
	From the above estimates, we can further bound $\partial_r\gamma_{k\ell}(r_1)-\partial_r\gamma_{k\ell}(r_2)$ as follows:
	\begin{align}
		|\partial_r \gamma_{k\ell}(r_1) - \partial_r \gamma_{k\ell}(r_2)| 
		&\le C\left[\left(\frac{1}{r_2}\right)\left|\gamma_{k\ell}(r_1)-\gamma_{k\ell}(r_2)\right| + \left(\frac{ \mu_k^{1/2}}{r_1^2r_2^3}\right)\left|r_1^3-r_2^3\right|\right]\notag\\
		&\le C\left[\left|\gamma_{k\ell}(r_1)-\gamma_{k\ell}(r_2)\right| + \mu_k^{1/2}\left|r_1^3-r_2^3\right|\right], \qquad\quad (r_1,r_2\ge 1).\notag
	\end{align}
	Lastly, we apply Lemmas \ref{l:equiv_lem0} and  \ref{l:equiv_lem2} to obtain:
	\begin{align*}
		|\partial_r \gamma_{k\ell}(r_1) - \partial_r \gamma_{k\ell}(r_2)| &\le C(1+\mu_k)^{1/2}\left(|r_1-r_2| + |r_1^3-r_2^3|\right)
	\end{align*}
\end{proof}

\begin{LemmaApdx}\label{l:equiv_lem5}
	For any $u\in H^1(S^{n-1},\C)$, there is a constant $C\ge0$ and independent of $r$ such that
	$$\frac{1}{r^2}\|u\|_{H^1}^2 \le 2C\|u\|_{L^2}^2+2\|A(r)^{1/2}u\|_{L^2}^2.$$
	Furthermore, if both \textup{(H2)} and \textup{(H3)} hold, then for some constant $C>0$ independent of $r$,
	$$\|u\|_{L^2}^2\le C\|A(r)^{1/2}u\|_{L^2}^2.$$		
\end{LemmaApdx}

\begin{proof}
	By expanding the $H^1$ norm of $u$, adding and subtracting $r^2 V_\infty$, and applying both triangle and Young's inequalities we obtain
	\begin{align*}
		\frac{1}{r^2}\lp\|u\rp\|_{H^1}^2 
		&= \sum_{k,j}  \bigg\|\bigg(\frac{1+\mu_k}{r^2}\bigg)^{1/2}u_{kj}\bigg\|^2_{\C^d} 
		\le \sum_{k,j}\lp\|\left(\left(1 + \mu_k/r^2\right)^{1/2}\pm A_k(r)^{1/2}\right) u_{kj}\rp\|_{\C^d}^2\\
		&\le  \sum_{k,j}  \Bigg[\left\|\lp(\lp(1+\mu_k/r^2\rp)^{1/2}-A_k(r)^{1/2}\rp)u_{kj}\right\|_{\C^d}+\left\|A_k(r)^{1/2}u_{kj}\right\|_{\C^d}\Bigg]^2\\
		&\le 2\bigg[ \sum_{k,j} \lp\|(1+\mu_k/r^2)^{1/2}-A_k(r)^{1/2}\rp\|^2_{L(\C^d)}\lp\|u_{kj}\rp\|_{\C^d}^2 \bigg]+ 2\lp\|A(r)^{1/2}u\rp\|_{L^2}^2
	\end{align*}
	We proceed by showing that the matrix norm of $(1+\mu_k/r^2)^{1/2}-A_k(r)^{1/2}$ is bounded above by a finite constant $C\ge 0$. Rewriting the matrix norm, we get
	\begin{align*}
		\left\|(1+\mu_k/r^2)^{1/2}-A_k(r)^{1/2}\right\| 
		& = \left\|\frac{1+\mu_k/r^2-V_\infty-\mu_k/r^2}{(1+\mu_k/r^2)^{1/2}+A_k(r)^{1/2}}\right\| 
		= \frac{\left\|1-V_\infty\right\|}{\left\|(1+\mu_k/r^2)^{1/2}+A_k(r)^{1/2}\right\|}. 
	\end{align*}
	Finally, since $\mu_k/r^2\ge0$ the denominator has a norm of at least $1$, and we have
	\begin{align*}
		\left\|(1+\mu_k/r^2)^{1/2}-A_k(r)^{1/2}\right\| 
		& \le \|1-V_\infty\| .
	\end{align*}
	Setting $C=\|1-V_\infty\|\ge 0$, the first statement follows. 
	
	To prove the second statement, notice that if (H2) holds, then we have an $r$-independent lower bound on $|\gamma_{k\ell}(r)|^2$  given by two cases:
	\begin{equation}\label{eq.eig_lowerbd}
	|\gamma_{k\ell}(r)|^2 = |\gamma_{k\ell}(r)^2| = \sqrt{\lp(\Re\lp[{\lambda_\ell}\rp]+\mu_k/r^2\rp)^2 + \Im\lp[{\lambda_\ell}\rp]^2} \ge \begin{cases} \big|\Im\lp[{\lambda_\ell}\rp]\big| & \text{ if } \Im\lp[{\lambda_\ell}\rp]\ne 0\\
		\,\big|\Re\lp[{\lambda_\ell}\rp]\big|	&\text { if } \Im\lp[{\lambda_\ell}\rp] = 0
	\end{cases}
	\end{equation}
	where we note that the inequality $\big|\gamma_{k\ell}(r)\big|^2\ge \big|\Re\lp[\lambda_\ell\rp]\big|$ holds by (H2) only in the case where $\Im\lp[\lambda_\ell\rp]=0$.

	If $d=1$, then \eqref{eq.eig_lowerbd} directly implies $\|A(r)^{1/2}u\|_{L^2}^2 \ge C \|u\|_{L^2}^2$ for a constant $C>0$ independent of $r$.
	
	If $d>1$, then \eqref{eq.eig_lowerbd} means that $|\gamma_{k\ell}(r)|$ is bounded away from zero uniformly in $k$, $\ell$, and $r$. Let $\Gamma>0$ be this uniform lower bound for all $|\gamma_{k\ell}(r)|$, and note that $1\le \frac{|\gamma_{k\ell}(r)|}{\Gamma}$. Using this bound, and the decomposition of the $L^2$ norm given in \eqref{eq.decomp_norms}, we have
	\begin{align*}
		\big\|u\big\|_{L^2}^2 &= \sum_{k,j}\|u_{kj}\|_{\C^d}^2 = \sum_{k,j}\|R R^{-1} u_{kj}\|_{\C^d}^2
		\le \sum_{k,j}\|R\|^2\big\|R^{-1}u_{kj}\big\|_{\C^d}^2
	\end{align*}
	where $R$ is the matrix whose columns are the eigenvectors $\{\rho_1,\ldots,\rho_d\}$ of $V_\infty$, guaranteed by (H3). Let $u_{kj}=\sum_{\ell=1}^d u_{kj}^{(\ell)}\rho_\ell$. Then $R^{-1}u_{kj}=(u_{kj}^{(1)},\ldots,u_{kj}^{(d)})\in\C^d$, and we have
	\[\|R^{-1}u_{kj}\|^2_{\C^d} = \sum_{\ell=1}^d \lp|u_{kj}^{(\ell)}\rp|^2.\]
	We proceed by bounding $\|u\|_{L^2}^2$ as follows:
	\begin{align*}
		\big\|u\big\|_{L^2}^2 &=   \sum_{k,j}\|R\|^2	\sum_{\ell=1}^d \lp|u_{kj}^{(\ell)}\rp|^2\\
		&\le 	 \sum_{k,j}\|R\|^2	\sum_{\ell=1}^d\frac{|\gamma_{k\ell}(r)|^2}{\Gamma^2}  \lp|u_{kj}^{(\ell)}\rp|^2\\
		&= \sum_{k,j}\frac{\|R\|^2}{\Gamma^2} \lp\| \begin{pmatrix}
			\gamma_{k1}(r) & &\\[-1em]
			&\!\!\!\!\ddots &\\[-1em]
			& & \!\!\!\!\gamma_{kd}(r)
		\end{pmatrix}R^{-1}u_{kj}\rp\|^2_{\C^d}
	\end{align*}
	
	Finally, recall that by the diagonalization of $V_\infty$, (H3),  we have that $R^{-1}A_{k}(r)^{1/2} = \diag(\gamma_{k\ell}(r))R^{-1}$. Therefore,
	\begin{align*}
		\big\|u\big\|_{L^2}^2 &\le \sum_{k,j}\frac{\|R\|^2}{\Gamma^2} \lp\|R^{-1}A_k(r)^{1/2}u_{kj}\rp\|^2_{\C^d} \le \frac{\|R\|^2\|R^{-1}\|^2}{\Gamma^2} \lp\|A(r)^{1/2}u\rp\|_{\C^d}^2.
	\end{align*}
\end{proof}

\begin{LemmaApdx}\label{l.real_lower_bd}
	The real part of $\gamma_{k\ell}(r)$ is bounded below by $\mathrm{Re}\lp[\lambda_\ell^{1/2}\rp]$.
	\[\mathrm{Re}\Big[\gamma_{k\ell}(r)\Big] = \mathrm{Re}\!\lp[\lp(\lambda_\ell+\mu_k/r^2\rp)^{1/2}\rp] \ge \mathrm{Re}\!\lp[\lambda_\ell^{1/2}\rp]\]
\end{LemmaApdx}

\begin{proof}
	If $\lambda_\ell\in\R$, then $\lambda_\ell>0$ by hypothesis (H2) and the result follows directly.
	
	If $\Im\lp[\lambda_\ell\rp] \ne 0$, then standard calculation of a complex principal square root gives 
		\[\mathrm{Re}\!\lp[\lp(\lambda_\ell+\mu_k/r^2\rp)^{1/2}\rp]^2 = \frac{1}{2}\Big(\lp|\lambda_\ell + \mu_k/r^2\rp|+\mathrm{Re}\!\lp[\lambda_\ell+\mu_k/r^2\rp]\Big) 
		 = \frac{g_\ell(\mu_k/r^2)}{2}\]
	where we define 
	$$g_\ell(\xi):= \left[\lp(\Re\lp[\lambda_\ell\rp] + \xi\rp)^2+ \Im\lp[\lambda_\ell\rp]^2\right]^{1/2} + (\Re\lp[\lambda_\ell\rp] + \xi).$$
	 A straightforward calculation shows that $g_\ell(\xi)$ is strictly increasing for all $\xi\ge 0$ with $\Im\lp[\lambda_\ell \rp]\ne0$.
 Therefore we have
	\begin{align*}
		\mathrm{Re}\!\lp[\lp(\lambda_\ell+\mu_k/r^2\rp)^{1/2}\rp]^2 =\frac{g_\ell(\mu_k/r^2)}{2}\ge\frac{g_\ell(0)}{2}= \frac{1}{2}\Big(\lp|\lambda_\ell\rp|+\Re\lp[\lambda_\ell\rp]\Big) = \mathrm{Re}\!\lp[\lambda_\ell^{1/2}\rp]^2. 
	\end{align*}
	
	The result follows since we are taking the principal square root.
\end{proof}

\section*{Acknowledgments}
MB and AHH gratefully acknowledge NSF support under award number DMS-2205434. RG gratefully acknowledges NSF support under Award Numbers DMS-2006887 and DMS-2307650.

\bibliographystyle{abbrvnat}
\bibliography{ses-exp-dich}

\end{document}